\documentclass[11pt]{article}
\usepackage[margin=1in]{geometry}
\usepackage[T1]{fontenc}
\usepackage[utf8]{inputenc}
\usepackage{lmodern}
\usepackage{microtype}
\DisableLigatures{encoding = T1, family = lmtt}
\usepackage{needspace}
\usepackage{amsmath,amssymb,amsthm,mathtools}
\usepackage{booktabs,array,longtable,tabularx}
\usepackage{enumitem}
\usepackage{xcolor}
\usepackage[hidelinks]{hyperref}
\usepackage{doi}
\usepackage{xurl}
\usepackage[backend=biber,citestyle=numeric-comp,bibstyle=ieee,sorting=none,minbibnames=5,maxbibnames=8,giveninits=true]{biblatex}
\AtBeginBibliography{\footnotesize\setlength{\bibitemsep}{0.25\baselineskip}}
\DeclareFieldFormat{doi}{\href{https://doi.org/#1}{doi:\nolinkurl{#1}}}
\hypersetup{
  pdftitle={Multivariable Congruence Jumping: Periodic States, Integral Realizations, and Descent},
  pdfauthor={Max A. Alekseyev and Dmitry I. Khomovsky},
  pdfsubject={Multivariable polynomial divisibility states, exponent cycles, integral realizations, and height descent},
  pdfkeywords={congruence jumping, Vieta mutations, Markoff-Hurwitz, integral realizations, Diophantine descent}
}
\newtheorem{theorem}{Theorem}[section]
\newtheorem{proposition}[theorem]{Proposition}
\newtheorem{lemma}[theorem]{Lemma}
\newtheorem{corollary}[theorem]{Corollary}
\newtheorem{example}[theorem]{Example}

\theoremstyle{definition}
\newtheorem{definition}[theorem]{Definition}
\theoremstyle{remark}
\newtheorem{remark}[theorem]{Remark}

\newcommand{\defeq}{:=}
\newcommand{\Z}{\mathbb Z}
\newcommand{\Q}{\mathbb Q}
\newcommand{\Ftwo}{\mathbb F_2}
\newcommand{\bx}{\mathbf x}
\newcommand{\widehatx}[1]{\bx_{\widehat{#1}}}

\title{Multivariable Congruence Jumping:\\Periodic States, Integral Realizations, and Descent}
\author{Max A. Alekseyev\thanks{The George Washington University, Washington, DC, USA. Corresponding author. Email: \href{mailto:maxal@gwu.edu}{maxal@gwu.edu}}
\and
Dmitry I. Khomovsky\thanks{Independent researcher. Email: \href{mailto:dmkhomovsky@gmail.com}{dmkhomovsky@gmail.com}}}
\date{}

\begin{document}
\maketitle

\begin{abstract}
Congruence jumping replaces coordinates by integer quotients and updates
polynomial divisibility relations. We extend this procedure to several
variables and retain integer coefficients by recording their common
factor separately. For three-variable systems of divisibilities with
signed binomial numerators, we classify all positive integer exponent
tuples of period at most two under three successive coordinate jumps:
seven fixed tuples and twenty-two cycles of exact period two. No odd
period greater than one occurs. At the chosen representatives of the twenty-two exponent two-cycles,
there are eighty-four returning sign assignments. All have realizations
with nonzero pairwise-coprime integer coordinates, and exactly
seventy-four have such realizations with every coordinate positive
throughout the six jumps. Explicit constructions and congruence arguments
establish this classification. We also study jumps that decrease a
specified height, distinguishing points that admit no immediate decrease
from minima of entire connected components of the jump graph. For the
four-variable family $Q_\lambda(x,y,z,w)=kxyzw$, where $Q_\lambda$ is the
sum of the four squares plus $\lambda$ times the six pairwise products,
finite bounds give complete classifications for $\lambda=1,2,3$.
For $\lambda=1$, positive integer solutions exist exactly for $k=2,8,10$.
These results distinguish returns of the polynomial data from integrality
and positivity of the coordinate sequences.
\end{abstract}



\section{Introduction}

Vieta jumping replaces one root of a quadratic equation by the other while
preserving the equation. In congruence jumping the basic requirement is
instead a collection of polynomial divisibilities, and replacing a
coordinate by a complementary quotient generally changes the polynomial
relations. The arithmetic of separated two-variable systems is developed
in~\cite{AlekseyevKhomovsky2026Arithmetic}; their coefficient dynamics is
studied in~\cite{AlekseyevKhomovsky2026Dynamics}. The present paper treats
systems
\begin{equation}\label{eq:intro-nvar}
  x_i\mid A_i(x_1,\ldots,\widehat{x_i},\ldots,x_n),
  \qquad 1\le i\le n.
\end{equation}
The ordered list $(A_1,\ldots,A_n)$ is the \emph{polynomial state}.
A \emph{realization} is an integer point with all coordinates nonzero
that satisfies these divisibilities. To update the relation for $x_i$
after a jump in $x_j$, with $i\ne j$, we set $x_i=0$ in $A_j$ and use the
resulting polynomial in a reciprocal substitution. The degrees, as well
as the coefficients, can therefore change.

The central question is how periodicity of the polynomial data relates
to actual integer solutions. We study this question for three-variable
systems of divisibilities with signed binomial numerators
\begin{equation}\label{eq:intro-binomial}
  A(y,z)=y^a+\alpha z^b,\qquad
  B(x,z)=x^c+\beta z^d,\qquad
  C(x,y)=x^e+\gamma y^f,
\end{equation}
where the exponents are positive and $\alpha,\beta,\gamma\in\{\pm1\}$.
We call the ordered triple $(A,B,C)$ a \emph{signed binomial state}:
\emph{binomial} refers to the two monomial terms in each polynomial,
and \emph{signed} records the three choices of sign. The six-tuple
$(a,b,c,d,e,f)$ is its \emph{exponent state}. A return of the exponents
need not return the signs, and a returning signed polynomial state need
not have a positive pairwise-coprime realization. The classifications below
distinguish these three questions.

\paragraph{Transport and exponent dynamics.}
Theorem~\ref{thm:face-transport} gives the reciprocal substitution that
updates the polynomials after one coordinate jump. The weighted version,
Theorem~\ref{thm:weighted-primitive-transport}, records the gcd of a
polynomial's coefficients separately from the polynomial obtained by
dividing out that gcd. For nonzero transformed polynomials, this requires
only exact division and coefficient gcds, with no integer factorization
or normalization to leading coefficient $1$.
At pairwise-coprime integer realizations of binomial states, the additional
monomial factors can be cancelled from the divisibility relations.
When all updated exponents remain positive, writing them in cyclic order
as $(q_1,\ldots,q_6)=(a,b,d,c,e,f)$ gives updates of the form
$q_i\mapsto q_{i-1}q_{i+1}-q_i$, with indices read modulo $6$.
These are Vieta mutations. On the symmetric restriction of
Section~\ref{sec:markoff-slice}, they reduce to Markoff--Hurwitz
transformations~\cite{Baragar1994,Baragar1998}. The state with all six
exponents equal to $2$ is fixed, whereas equal initial exponents
$r\ge3$ lead to strict growth.

One jump in $x$, followed by one in $y$ and one in $z$, will be called an
$xyz$ block. Theorem~\ref{thm:period2-param} reduces period two of the
resulting map to a finite unit-fraction problem. There are exactly seven
fixed exponent states and twenty-two cycles of exact period two
(Corollary~\ref{cor:22-cycles}). A quantity that changes sign after each
block, together with a positive recurrence, excludes every odd period
greater than one (Theorem~\ref{thm:no-odd-period}). The classification
concerns periods one and two; higher even periods remain open.

\paragraph{Signs and arithmetic realizations.}
After encoding signs as $0$ and $1$, their updates are linear over the
two-element field $\Ftwo$. For each of the twenty-two exponent cycles,
choose the representative listed in the first exponent column of
\ref{app:cycles}. Theorem~\ref{thm:84-signs} gives exactly eighty-four
returning sign assignments at these representatives. The corresponding
polynomial states all have pairwise-coprime nonzero integral
six-jump realizations. Exactly seventy-four have positive pairwise-coprime realizations
(Theorem~\ref{thm:positive-classification}); the other ten are excluded
by congruences at the point reached after three jumps and inequalities
forced by the negative signs. For the positive cases, the existence
proof reconstructs six-jump sequences by solving the quotient equations
backwards from explicitly chosen data. These constructions establish one
six-jump return; they do not by themselves establish positivity under
arbitrarily many repetitions.

For a pairwise-coprime point, the three relations can also be written as
$F_{A,B,C}(x,y,z)=kxyz$, with $F_{A,B,C}=yzA+xzB+xyC$.
Example~\ref{ex:k-change} shows that the same polynomial can return after
six jumps with a different integer value of $k$. Thus a return of the polynomial
data need not preserve the value of $k$.

\paragraph{Roots and descent.}
Section~\ref{sec:reverse-descent} studies a graph whose vertices consist
of polynomial data together with an integer solution, and whose edges
are the permitted jumps. Heights compare the sorted absolute coordinate
values. A \emph{local root} admits no height-decreasing jump, but need
not minimize height in its connected component: changing the polynomial
data while keeping the coordinates fixed can permit a later decrease.
Proposition~\ref{prop:root-descent} also treats roots that remain local
roots after any sequence of height-preserving jumps. No periodicity of
the polynomial data is assumed.
For the four-variable family $Q_\lambda=kxyzw$, where $Q_\lambda$ is the
sum of the four squares plus $\lambda$ times all pairwise products,
Theorem~\ref{thm:root-strip} bounds roots with sorted positive coordinates
$x\le y\le z\le w$ by $\lambda+3\le kxy\le6\lambda+4$.
This gives a finite classification for each fixed positive integer
$\lambda$. For $\lambda=1$, positive solutions exist exactly for
$k\in\{2,8,10\}$; the root lists for $\lambda=2,3$ are also complete.

\paragraph{Relation to other work and scope.}
The polynomial substitutions are related to cluster and Laurent
mutations~\cite{FominZelevinsky2002,LamPylyavskyy2016}. However, our
cancellations depend on divisibility and coprimality at the integer
point. A Laurent-polynomial expression, which may contain negative
powers, does not by itself ensure that a quotient jump is integral. The cross-term family is compared with generalized
Hurwitz equations~\cite{FineEtAl2016Hurwitz}, cluster transformations
\cite{Lampe2016,BaoLi2025}, and higher-dimensional mutation theory
\cite{Cho2026} in Section~\ref{sec:reverse-descent}.

For an acyclic cluster seed with trivial coefficients, positive
integrality of the initial variables and of their one-step replacements
ensures positive integrality of every cluster variable
\cite[Section~6.3]{deSaintGermainHuangLu2026}.
Zhang~\cite[Theorem~B]{Zhang2025} determines when the associated
varieties have finitely many positive integral points. His rank-two examples include
$x\mid y+1$ and $y\mid(x^r+1)^s$ for positive integers $r,s$.
Each quotient jump fixes the polynomial pair $(t+1,(t^r+1)^s)$, and
\cite[Section~5.1 and Remark~5.1]{Zhang2025} gives a growing iteration
when $rs\ge4$. The map from positive integral cluster points to
$xyz=(x^r+1)^s+y$ is injective, with image satisfying $y=u^s$ for a
positive integer $u$; it is bijective when $s=1$.
Proposition~\ref{prop:cluster-periodic} compares the ordinary
exchange-matrix condition with periodic exponent states whose
intermediate exponents remain positive: it forces a fixed state, and
the matrices of all seven fixed states contain a directed cycle after
every mutation.

Backward construction from a chosen final integer point follows the
classical Mordell--Schinzel procedure, with integrality of the reverse
steps proved separately
\cite[equations~(50)--(53) and Lemma~14]{Schinzel2015}.
A positive-monomial version with explicit growth estimates is given in
\cite[Subsection~2.6, Proposition~2.19]{AlekseyevKhomovsky2026Arithmetic}.
Section~\ref{sec:face-transport} gives a sufficient reverse-step condition
in several variables, and Section~\ref{sec:positive-lifts} implements it
for finite three-variable realizations. Neither a backward chain of
arbitrary prescribed length nor the ability to continue jumps forwards
is sufficient, by itself, to produce infinitely many distinct solutions
of the original divisibility system.

Using generalized cluster algebras of rank two, de Saint Germain and
Zhang~\cite{deSaintGermainZhang2026} study positive integral points on
surfaces $uvt=G(u,v)$, where $G$ has nonnegative integer coefficients. When $G(0,0)>0$, their reciprocal transforms of the
one-variable polynomials $G(u,0)$ and $G(0,v)$ agree, up to constant
factors, with the two-variable specialization of
Theorem~\ref{thm:face-transport}. The two variables
replaced by quotient jumps are $u$ and $v$; $t$ records $G(u,v)/(uv)$.
Their finiteness, Zariski density, and orbit results give geometric context for
this specialization. The classifications here concern returns of the
six exponents and three signs in three-variable binomial states, together
with integral and positive realizations of those returns.

Section~\ref{sec:face-transport} establishes transport. Sections
\ref{sec:binomial}--\ref{sec:positive-lifts} develop exponent dynamics,
periods, signs, and integral realizations in that order.
Section~\ref{sec:reverse-descent} treats roots and the four-variable family.
The appendices give the complete finite tables, starting data for the
positive realizations, and the finite reductions used to verify
completeness.

\section{Reciprocal transport between polynomial states}\label{sec:face-transport}

Let $\bx=(x_1,\ldots,x_n)$, and write $\widehatx{i}$ for the tuple obtained by omitting $x_i$. For each $i$, let
\[
 A_i\in\Z[x_1,\ldots,\widehat{x_i},\ldots,x_n].
\]
An integer point $\bx\in(\Z\setminus\{0\})^n$ \emph{realizes} the polynomial state $(A_1,\ldots,A_n)$ if
\[
 x_i\mid A_i(\widehatx{i})
 \qquad(1\le i\le n).
\]

Fix $j$.  At a realizing point define the quotient jump
\begin{equation}\label{eq:j-jump}
 x_j'\defeq \frac{A_j(\widehatx{j})}{x_j}.
\end{equation}
For $i\ne j$ with $A_i\ne0$, write
\[
 d_{ij}\defeq\deg_{x_j}A_i.
\]
If $A_i=0$, set $d_{ij}\defeq0$. Define $S_{ij}$ by setting $x_i=0$ in $A_j$:
\begin{equation}\label{eq:face-specialization}
 S_{ij}\defeq A_j\big|_{x_i=0}.
\end{equation}
The polynomial $S_{ij}$ involves neither $x_i$ nor $x_j$.

\begin{theorem}\label{thm:face-transport}
For $i\ne j$, using $X_j$ for the new $j$th coordinate variable, define
\begin{equation}\label{eq:face-transform}
 A_i^{[j]}(\ldots,X_j,\ldots)
 \defeq
 X_j^{d_{ij}}
 A_i\!\left(\ldots,\frac{S_{ij}}{X_j},\ldots\right),
\end{equation}
and put $A_j^{[j]}\defeq A_j$.  Then every $A_i^{[j]}$ has integer coefficients.  If $\bx$ realizes $(A_1,\ldots,A_n)$ and $x_j'\ne0$ is given by \eqref{eq:j-jump}, then the jumped point
\[
 (x_1,\ldots,x_{j-1},x_j',x_{j+1},\ldots,x_n)
\]
realizes $(A_1^{[j]},\ldots,A_n^{[j]})$.
\end{theorem}

\begin{proof}
Fix $i\ne j$ and expand
\[
 A_i=\sum_{r=0}^{d_{ij}} a_r x_j^r,
\qquad
a_r\in\Z[x_1,\ldots,\widehat{x_i},\ldots,\widehat{x_j},\ldots,x_n].
\]
Then
\[
 A_i^{[j]}=\sum_{r=0}^{d_{ij}}a_rS_{ij}^rX_j^{d_{ij}-r},
\]
so polynomiality and integrality are immediate.

Since $x_jx_j'=A_j(\widehatx{j})$, reduction modulo $x_i$ gives
\[
 x_jx_j'\equiv S_{ij}\pmod{x_i}.
\]
Hence
\[
\begin{aligned}
 A_i^{[j]}(x_j')
 &=\sum_{r=0}^{d_{ij}}a_rS_{ij}^r(x_j')^{d_{ij}-r}\\
 &\equiv
 (x_j')^{d_{ij}}
 \sum_{r=0}^{d_{ij}}a_rx_j^r
 =
 (x_j')^{d_{ij}}A_i(\widehatx{i})
 \equiv0\pmod{x_i}.
\end{aligned}
\]
The $j$th divisibility is preserved because $x_jx_j'=A_j(\widehatx{j})$.
\end{proof}

\paragraph{Backward reconstruction.}
The congruence in the proof also gives a sufficient condition for reconstructing a predecessor. Suppose that the preceding polynomials $A_i$ and the jump index $j$ are known. Let $\bx'$ realize $(A_1^{[j]},\ldots,A_n^{[j]})$, set $x_i:=x_i'$ for $i\ne j$, and define
\[
 x_j\defeq\frac{A_j(\widehatx{j})}{x_j'}.
\]
This quotient is integral because $A_j^{[j]}=A_j$; assume it is nonzero.
The product relation $x_jx_j'=A_j(\widehatx{j})$ gives
\[
 A_i^{[j]}(\widehatx{i}')
 \equiv (x_j')^{d_{ij}}A_i(\widehatx{i})\pmod{x_i}
 \qquad(i\ne j),
\]
where $\widehatx{i}'$ denotes $\bx'$ with its $i$th coordinate omitted.
If $\gcd(x_j',x_i)=1$ for every $i\ne j$, the displayed congruence
implies $x_i\mid A_i(\widehatx{i})$; the $j$th divisibility follows from
the product identity. Hence $\bx$ realizes the preceding state.
This is reconstruction along a specified sequence, not an assertion
that applying the forward coefficient transform twice restores the
original state.  For a sequence of reverse steps, the nonzero quotients
and the required cancellations must remain valid at every stage.
After a normalization, its retained weights or additional cancellation
conditions must also be respected.  Thus backward construction is
available in the multivariable setting, but an arbitrary solution of a
chosen final state need not permit a backward jump.

\begin{remark}\label{rem:face-hierarchy}
The parameter $S_{ij}$ in \eqref{eq:face-transform} is obtained by setting $x_i=0$ in $A_j$, so it depends on at most $n-2$ variables. For $n=2$ it is a scalar, reproducing the reciprocal-polynomial transport of~\cite{AlekseyevKhomovsky2026Arithmetic}. For $n=3$ it is generally a one-variable polynomial. Its nonconstant powers in \eqref{eq:face-transform} explain why transport can change the degrees in the other variables.
\end{remark}

\Needspace{6\baselineskip}
A single polynomial equation supplies the divisibilities used in transport.
\begin{corollary}\label{cor:F-source}
Let $F\in\Z[x_1,\ldots,x_n]$ and $k\in\Z$. If
\[
 F(x_1,\ldots,x_n)=k\,x_1\cdots x_n
\]
at $\bx\in(\Z\setminus\{0\})^n$, then
\[
 x_i\mid F\big|_{x_i=0}
 \qquad(1\le i\le n).
\]
\end{corollary}

\begin{proof}
Reduce the equation modulo $x_i$.
\end{proof}

The converse holds for some integer $k$ when the coordinates are pairwise coprime. The displayed divisibilities imply $x_i\mid F(\bx)$ for every $i$, since $F(\bx)\equiv F|_{x_i=0}(\bx)\pmod{x_i}$. Pairwise coprimality then gives $x_1\cdots x_n\mid F(\bx)$.  Without this hypothesis the separate divisibilities need
not imply divisibility by the product.  For example, at $(2,2,2)$ the
polynomial $F(x,y,z):=x^3+y^3+z^3+2$ takes the value $26$, which is
divisible by each coordinate but not by $xyz=8$.  The general state
transport concerns the separate relations; it does not require this
converse or pairwise coprimality.

\subsection{Transport with integer weights}\label{subsec:weighted-primitive}

Theorem~\ref{thm:face-transport} transports integer polynomials without any normalization. For a nonzero integer polynomial $P$, its content $\operatorname{cont}(P)$ is the positive gcd of its coefficients; $P$ is primitive if this gcd is $1$. We separate each polynomial into its content and primitive part, retaining the content as a positive integer weight instead of dividing by a prescribed leading or constant coefficient.

\begin{definition}\label{def:weighted-primitive-state}
For each $i$, let $\lambda_i\in\Z_{>0}$ and let
\[
 A_i\in\Z[x_1,\ldots,\widehat{x_i},\ldots,x_n]
\]
be a nonzero primitive polynomial.  An integer point $\bx\in(\Z\setminus\{0\})^n$ \emph{realizes the weighted state}
\[
 ((\lambda_1,A_1),\ldots,(\lambda_n,A_n))
\]
if
\begin{equation}\label{eq:weighted-state-divisibility}
 x_i\mid \lambda_i A_i(\widehatx{i})
 \qquad(1\le i\le n).
\end{equation}
The scalar $\lambda_i$ is called the \emph{weight} of the $i$th relation.  An overall sign may always be absorbed into $A_i$, so taking the weights positive loses no information.
\end{definition}

Let $\bx$ realize the weighted state. Fix a jump in the $j$th coordinate
and assume its quotient is nonzero:
\begin{equation}\label{eq:weighted-jump}
 x_j'\defeq\frac{\lambda_jA_j(\widehatx{j})}{x_j}
 \in\Z\setminus\{0\}.
\end{equation}
For $i\ne j$ put
\begin{equation}\label{eq:weighted-face-parameter}
 C_{ij}\defeq
 \lambda_jA_j\big|_{x_i=0},
 \qquad
 d_{ij}\defeq\deg_{x_j}A_i,
\end{equation}
and form the transformed polynomial before extracting its content:
\begin{equation}\label{eq:weighted-raw-transform}
 Q_i^{[j]}(\ldots,X_j,\ldots)
 \defeq
 \lambda_iX_j^{d_{ij}}
 A_i\!\left(\ldots,\frac{C_{ij}}{X_j},\ldots\right).
\end{equation}

\begin{theorem}\label{thm:weighted-primitive-transport}
Every $Q_i^{[j]}$ in~\eqref{eq:weighted-raw-transform} has integer coefficients, and the jumped point satisfies
\begin{equation}\label{eq:weighted-raw-divisibility}
 x_i\mid Q_i^{[j]}(\ldots,x_j',\ldots)
 \qquad(i\ne j).
\end{equation}
The $j$th relation remains
\[
 x_j'\mid\lambda_jA_j(\widehatx{j}).
\]
Assume $Q_i^{[j]}\ne0$ for every $i\ne j$, and define
\begin{equation}\label{eq:weighted-content-extraction}
 g_i\defeq\operatorname{cont}(Q_i^{[j]})>0,
 \qquad
 A_i'\defeq Q_i^{[j]}/g_i,
 \qquad
 \lambda_i'\defeq g_i.
\end{equation}
Then $A_i'$ is primitive and the jumped point realizes the new weighted state obtained by replacing $(\lambda_i,A_i)$ with $(\lambda_i',A_i')$ for $i\ne j$ and leaving $(\lambda_j,A_j)$ unchanged.  When every $Q_i^{[j]}$ is nonzero, this normalization uses exact integer division and gcds of coefficients, without integer factorization.
\end{theorem}

\begin{proof}
Fix $i\ne j$ and expand
\[
 A_i=\sum_{r=0}^{d_{ij}}a_rx_j^r.
\]
Then
\[
 Q_i^{[j]}
 =\lambda_i\sum_{r=0}^{d_{ij}}
 a_rC_{ij}^rX_j^{d_{ij}-r},
\]
so it is an integer polynomial.  From~\eqref{eq:weighted-jump}, reduction modulo $x_i$ gives
\[
 x_jx_j'
 =\lambda_jA_j(\widehatx{j})
 \equiv C_{ij}\pmod{x_i}.
\]
Consequently
\[
 Q_i^{[j]}(x_j')
 \equiv
 \lambda_i(x_j')^{d_{ij}}A_i(\widehatx{i})
 \equiv0\pmod{x_i},
\]
which proves~\eqref{eq:weighted-raw-divisibility}.  The $j$th divisibility follows from
$\lambda_jA_j(\widehatx{j})=x_jx_j'$.  Finally,~\eqref{eq:weighted-content-extraction} changes only the presentation of the same integer polynomial $Q_i^{[j]}=\lambda_i'A_i'$, so it preserves every divisibility relation.  Computing a content is a gcd computation on coefficients and does not require their prime factorization.
\end{proof}

The weighted theorem constructs a successor state when the jump quotient and the transformed polynomials are nonzero. Keeping the content as a weight preserves the same polynomial relation, but does not make the transformation a bijection on all integral solutions. A backward construction must also verify the preceding divisibilities and recover nonzero integer coordinates. Even if this can be done for arbitrarily many steps, a separate argument is needed to show that the recovered solutions are distinct. The growth bounds and counterexamples in \cite[Subsection~2.6]{AlekseyevKhomovsky2026Arithmetic} distinguish these questions.

\begin{remark}\label{rem:weighted-zero-relation}
Since $A_i\ne0$, the transformed polynomial $Q_i^{[j]}$ can be identically zero only when $C_{ij}=0$. In that case,
\[
 Q_i^{[j]}=\lambda_iX_j^{d_{ij}}A_i\big|_{x_j=0},
\]
so it is identically zero exactly when $A_i|_{x_j=0}=0$. The corresponding divisibility is then automatic, and no nonzero primitive polynomial can be extracted. The unnormalized transformation remains valid; the content normalization in Theorem~\ref{thm:weighted-primitive-transport} applies only to nonzero transformed polynomials.
\end{remark}

For two-variable states $(A(y),B(x))$, the specialization parameters are scalars. The corresponding weighted transport and return after four alternating jumps are proved in the companion Dynamics paper~\cite[Section~2]{AlekseyevKhomovsky2026Dynamics}. The return is up to multiplication and rescaling of the variable: $Q(T)=uP(vT)$ with $u,v\in\Q^\times$. It assumes nonzero constant terms in both initial polynomials, so each reciprocal update preserves degree and keeps the leading and constant coefficients nonzero. Theorem~\ref{thm:weighted-primitive-transport} itself does not require nonzero constant terms. In several variables the parameters are generally nonconstant, so the two-variable result does not imply a return of the multivariable polynomial state.

\section{Three-variable signed binomial states}\label{sec:binomial}

We now specialize to the signed binomial states introduced in
\eqref{eq:intro-binomial}:
\begin{equation}\label{eq:ABC}
 A=y^a+\alpha z^b,\qquad
 B=x^c+\beta z^d,\qquad
 C=x^e+\gamma y^f,
\end{equation}
where
\[
 a,b,c,d,e,f\in\Z_{>0},
 \qquad
 \alpha,\beta,\gamma\in\{\pm1\}.
\]
An integer triple $(x,y,z)$ with all coordinates nonzero is \emph{pairwise coprime} if
\[
 \gcd(x,y)=\gcd(x,z)=\gcd(y,z)=1.
\]

\begin{lemma}\label{lem:coprime-preserved}
Suppose $(x,y,z)$ is pairwise coprime and
\[
 x'=\frac{y^a+\alpha z^b}{x}
\]
is a nonzero integer.  Then $(x',y,z)$ is pairwise coprime.  The analogous statements hold for $y$- and $z$-jumps.
\end{lemma}

\begin{proof}
If a prime $p$ divided both $x'$ and $y$, then
\[
 0\equiv xx'=y^a+\alpha z^b\equiv\alpha z^b\pmod p,
\]
so $p\mid z$, contradicting $\gcd(y,z)=1$.  The proof for $\gcd(x',z)$ is identical.
\end{proof}

The lemma allows common monomial factors to be cancelled in the transported divisibilities, because they are coprime to the corresponding coordinate.  For example, under an $x$-jump Theorem~\ref{thm:face-transport} gives
\[
 X^cB\!\left(\frac{\alpha z^b}{X},z\right)
 =
 \alpha^cz^{bc}+\beta z^dX^c.
\]
If $bc>d$, the factor $\beta z^d$ is coprime to $y$, so it can be cancelled from the congruence modulo $y$, giving
\[
 B^{(x)}(X,z)=X^c+\beta\alpha^c z^{bc-d}.
\]
Similarly, if $ae>f$,
\[
 C^{(x)}(X,y)=X^e+\gamma y^{ae-f}.
\]

Define
\begin{equation}\label{eq:first-differences}
 d_1\defeq bc-d,\qquad f_1\defeq ae-f,
\end{equation}
then
\begin{equation}\label{eq:second-differences}
 b_2\defeq ad_1-b,\qquad e_2\defeq cf_1-e,
\end{equation}
and finally
\begin{equation}\label{eq:third-differences}
 a_3\defeq f_1b_2-a,\qquad c_3\defeq d_1e_2-c.
\end{equation}

\begin{definition}\label{def:positive-chamber}
A signed binomial state \eqref{eq:ABC} lies in the positive exponent region of one complete $xyz$ block if
\[
 d_1,f_1,b_2,e_2,a_3,c_3>0.
\]
\end{definition}

\begin{proposition}\label{prop:block-mutation}
In the positive exponent region, a complete $x$-, then $y$-, then $z$-jump sends
\[
 p=(a,b,c,d,e,f)
\]
to
\begin{equation}\label{eq:T-p}
 T(p)=(a_3,b_2,c_3,d_1,e_2,f_1).
\end{equation}
The signs evolve by
\begin{equation}\label{eq:sign-block}
 \beta'=\beta\alpha^c,\qquad
 \gamma'=\gamma,\qquad
 \alpha'=\alpha(\beta\alpha^c)^a\gamma^{b_2}.
\end{equation}
\end{proposition}

\begin{proof}
The $x$-jump gives \eqref{eq:first-differences} and replaces $\beta$ by $\beta\alpha^c$.  The $y$-jump gives \eqref{eq:second-differences} and multiplies $\alpha$ by $(\beta\alpha^c)^a$.  The $z$-jump gives \eqref{eq:third-differences} and multiplies the current $\alpha$ by $\gamma^{b_2}$.  The sign of $\gamma$ is unchanged throughout.
\end{proof}

Reorder the exponents as
\begin{equation}\label{eq:q-order}
 q=(q_1,\ldots,q_6)\defeq(a,b,d,c,e,f).
\end{equation}
Each jump updates a pair of entries three positions apart in this cyclic order, which we call an opposite pair. The $x$-jump updates $(q_3,q_6)$, the $y$-jump $(q_2,q_5)$, and the $z$-jump $(q_1,q_4)$, always by
\begin{equation}\label{eq:cyclic-vieta}
 q_i\longmapsto q_{i-1}q_{i+1}-q_i
 \qquad (i\bmod6).
\end{equation}

\begin{remark}\label{rem:laurent-chamber}
If one of the differences in \eqref{eq:first-differences}--\eqref{eq:third-differences} is negative, the corresponding normalized expression contains a negative power. The unnormalized transform remains polynomial; factoring a different common monomial instead produces a polynomial outside the form \eqref{eq:ABC}. We restrict to the positive exponent region, where every intermediate state retains that form. Periodicity of the formal exponent updates alone does not guarantee an orbit of the form \eqref{eq:ABC} or an integral realization.
\end{remark}

\section{A Markoff--Hurwitz invariant and exponent growth}\label{sec:markoff-slice}

Consider exponent states of the form
\begin{equation}\label{eq:symmetric-p}
 (a,b,c,d,e,f)=(u,v,u,w,v,w).
\end{equation}
This form is preserved by an $xyz$ block. Define the successive exponent updates by
\begin{equation}\label{eq:uvw-mutations}
 w'=uv-w,\qquad
 v'=uw'-v,\qquad
 u'=v'w'-u,
\end{equation}
the next state is again of the form \eqref{eq:symmetric-p}, with $(u,v,w)$ replaced by $(u',v',w')$.

\begin{proposition}\label{prop:markoff-invariant}
Each of the three Vieta reflections in \eqref{eq:uvw-mutations} preserves
\begin{equation}\label{eq:H}
 H(u,v,w)\defeq u^2+v^2+w^2-uvw.
\end{equation}
\end{proposition}

\begin{proof}
For fixed $u,v$, the two roots in the $w$-coordinate of the equation $H(u,v,w)=h$ have sum $uv$.  Thus $w\mapsto uv-w$ exchanges the roots.  The other two reflections are identical after permutation.
\end{proof}

These exponent updates have the same Vieta form as the Markoff--Hurwitz transformations~\cite{Baragar1994,Baragar1998}. Here the transformations act on $u,v,w$, which are exponents in the original polynomials.

\begin{proposition}\label{prop:diagonal-growth}
The state $(u,v,w)=(2,2,2)$ is fixed. For every $r\ge3$, the state $(r,r,r)$ enters a strictly growing orbit.  More precisely, after one block the state satisfies
\[
 u_1>v_1>w_1>r,
\]
and whenever $u>v>w\ge3$, the next block satisfies
\[
 u'>v'>w'>u.
\]
\end{proposition}

\begin{proof}
For $r=2$, every mutation gives $2\cdot2-2=2$.  For $r\ge3$,
\[
 w_1=r^2-r>r,\qquad
 v_1=rw_1-r>w_1,\qquad
 u_1=v_1w_1-r>v_1.
\]
Now suppose $u>v>w\ge3$.  Then
\[
 w'=uv-w>u,
\]
because $u(v-1)>w$.  Next
\[
 v'=uw'-v>w',
\]
and finally $u'=v'w'-u>v'$.  The inequalities are strict.
\end{proof}

\begin{example}\label{ex:cubic}
Consider
\begin{equation}\label{eq:cubic-Fermat}
 x^3+y^3-z^3=kxyz.
\end{equation}
A solution with $\gcd(x,y,z)=1$ is automatically pairwise coprime: if a prime divided two variables, \eqref{eq:cubic-Fermat} would force it to divide the third.

The associated divisibility state is
\[
 A=y^3-z^3,\qquad
 B=x^3-z^3,\qquad
 C=x^3+y^3.
\]
We continue to write $x,y,z$ for the coordinate variables after each jump. After an $x$-jump and cancellation of common monomial factors,
\[
 B\longmapsto x^3+z^6,\qquad
 C\longmapsto x^3+y^6.
\]
The subsequent $y$- and $z$-jumps complete the following sequence of polynomial states:
\[
\begin{aligned}
&(y^3-z^3,\ x^3-z^3,\ x^3+y^3)\\
&\longrightarrow
(y^3-z^3,\ x^3+z^6,\ x^3+y^6)\\
&\longrightarrow
(y^3-z^{15},\ x^3+z^6,\ x^{15}+y^6)\\
&\longrightarrow
(y^{87}-z^{15},\ x^{87}+z^6,\ x^{15}+y^6).
\end{aligned}
\]
Thus $H(3,3,3)=0$, and the exponents increase strictly under repeated $xyz$ blocks by Proposition~\ref{prop:diagonal-growth}.

After a quotient jump, the cubic expression in \eqref{eq:cubic-Fermat} need not be an integer multiple of the new coordinate product. For instance, $(x,y,z)=(1,2,3)$ has $k=-3$ and $x'=-19$. Since $y^3-z^3=xx'$, the same cubic expression at $(x',2,3)$ would require
\[
 k'=\frac{x'^2+x}{yz}=\frac{181}{3}.
\]
This is why periodicity of the full polynomial state, rather than integrality of one quotient coordinate, is the relevant return condition.
\end{example}

\section{Periodic exponent states}\label{sec:periodic-exponents}

Throughout this section, exponent states have positive integer entries,
and every intermediate exponent state in an orbit under consideration
is required to remain positive. We classify the fixed points and cycles
of exact period two of $T$, and then exclude odd periods greater than one.

\subsection{Fixed points}

\begin{theorem}\label{thm:fixed-exponents}
There are exactly seven positive integer fixed points of $T$.  In cyclic order \eqref{eq:q-order} they are
\[
 (2,2,2,2,2,2)
\]
and the six cyclic rotations of
\[
 (1,1,2,4,4,2).
\]
\end{theorem}

\begin{proof}
If $Tq=q$, then the first mutated pair $(q_3,q_6)$ must already be fixed, because these coordinates are not touched again in the block.  The same applies successively to $(q_2,q_5)$ and $(q_1,q_4)$.  Hence
\begin{equation}\label{eq:fixed-local}
 2q_i=q_{i-1}q_{i+1}
 \qquad(i\bmod6).
\end{equation}
Multiplying all six equations gives
\[
 \prod_{i=1}^6q_i=64.
\]
Thus every $q_i$ is a power of $2$, say $q_i=2^{r_i}$ with $r_i\ge0$.  Equation \eqref{eq:fixed-local} becomes
\[
 r_{i-1}+r_{i+1}=r_i+1.
\]
With $s_i=r_i-1$ this is
\[
 s_{i-1}+s_{i+1}=s_i.
\]
Writing $s_1=u$, $s_2=v$ gives
\[
 (s_1,\ldots,s_6)=(u,v,v-u,-u,-v,u-v).
\]
The condition $r_i=s_i+1\ge0$ leaves the seven possibilities
\[
 (u,v)\in\{(0,0),(-1,-1),(-1,0),(0,1),(1,1),(1,0),(0,-1)\},
\]
which yield the displayed states.
\end{proof}

\subsection{The period-two parametrization}

\begin{theorem}\label{thm:period2-param}
Let $q=(q_1,\ldots,q_6)$ have positive integer entries.  Then $T^2q=q$ if and only if there exist positive integers
\[
 g,h,m,n,U,V,
 \qquad \gcd(m,n)=1,
\]
such that
\begin{equation}\label{eq:core-div}
 n\mid m+U,\qquad m\mid n+V,
\end{equation}
\begin{equation}\label{eq:unit-fraction}
 \left(\frac1m+\frac1U\right)
 \left(\frac1n+\frac1V\right)=gh,
\end{equation}
and
\begin{equation}\label{eq:q-param}
 q=
 \left(
 hm,\ \frac{m+U}{n},\ gm,\ gn,\ \frac{n+V}{m},\ hn
 \right).
\end{equation}
Applying $T$ once interchanges $(m,n)$ and $(U,V)$ while leaving $g,h$ fixed.
\end{theorem}

\begin{proof}
Put $r=Tq$ and assume $Tr=q$.  Comparing the forward and reverse mutation equations gives
\[
 \frac{r_1}{q_1}=\frac{r_3}{q_3},
 \qquad
 \frac{r_5}{q_5}=\frac{q_1}{r_1},
 \qquad
 \frac{r_4}{q_4}=\frac{r_6}{q_6},
 \qquad
 \frac{r_2}{q_2}=\frac{q_4}{r_4}.
\]
In particular
\[
 q_1q_4=q_3q_6.
\]
Since $q_3/q_4=q_1/q_6$, write their common value as $m/n$ in lowest terms.  Then
\[
 q_3=gm,\quad q_4=gn,\quad q_1=hm,\quad q_6=hn,
 \qquad \gcd(m,n)=1.
\]
The mutation equations then give
\[
 U\defeq q_2n-m>0,\qquad
 V\defeq mq_5-n>0,
\]
and hence
\[
 q_2=\frac{m+U}{n},\qquad
 q_5=\frac{n+V}{m}.
\]
The remaining mutation equations reduce to
\[
 q_2q_5=ghUV,
\]
which is equivalent to \eqref{eq:unit-fraction}.  Reversing these steps proves the converse, and direct substitution shows that one block interchanges $(m,n)$ with $(U,V)$.
\end{proof}

\begin{corollary}\label{cor:22-cycles}
There are exactly $51$ positive exponent states with period dividing two: the seven fixed states of Theorem~\ref{thm:fixed-exponents} and $44$ further states forming exactly $22$ cycles of exact period two.
\end{corollary}

\begin{proof}
The left side of \eqref{eq:unit-fraction} is at most $4$, so $G\defeq gh\in\{1,2,3,4\}$.  The finite enumeration in \ref{app:cores}, subject to the divisibilities in \eqref{eq:core-div}, gives respectively
\[
 20,\ 10,\ 4,\ 1
\]
ordered quadruples $(m,n,U,V)$ for $G=1,2,3,4$.  Allowing the ordered factorizations $G=gh$ gives
\[
 20+2\cdot10+2\cdot4+3\cdot1=51
\]
states.  Exactly seven have $(m,n)=(U,V)$ and are fixed; the remaining $44$ form $22$ pairs exchanged by $T$.
\end{proof}

\subsection{An odd-period obstruction}

For cyclic exponent coordinates $q_1,\ldots,q_6$ put
\begin{equation}\label{eq:ABXYZ}
 \mathcal A=q_1q_3q_5,\qquad
 \mathcal B=q_2q_4q_6,
\end{equation}
and
\[
 X=q_1q_4,\qquad
 Y=q_2q_5,\qquad
 Z=q_3q_6.
\]

\begin{proposition}\label{prop:AB-invariants}
Under one mutation of an opposite pair, $\mathcal A-\mathcal B$ changes sign, while
\begin{equation}\label{eq:I-invariant}
 \mathcal I\defeq
 \mathcal A+\mathcal B-2(X+Y+Z)
\end{equation}
is invariant.  Consequently
\[
 (\mathcal A-\mathcal B)(Tq)=-(\mathcal A-\mathcal B)(q),
 \qquad
 \mathcal I(Tq)=\mathcal I(q).
\]
\end{proposition}

\begin{proof}
For example, mutating $(q_3,q_6)$ gives
\[
 q_3'=q_2q_4-q_3,\qquad
 q_6'=q_5q_1-q_6.
\]
Since $X$ and $Y$ are unchanged,
\[
 \mathcal A'=XY-\mathcal A,\qquad
 \mathcal B'=XY-\mathcal B,
\]
so the difference changes sign.  Also
\[
 Z'=XY-\mathcal A-\mathcal B+Z,
\]
which gives \eqref{eq:I-invariant}.  A complete block has three opposite-pair mutations.
\end{proof}

The exclusion of odd periods uses the following recurrence lemma.
\begin{lemma}\label{lem:sqrt-recurrence}
Let $(u_j)_{j\in\Z}$ be a positive periodic sequence satisfying
\begin{equation}\label{eq:sqrt-recurrence}
 u_{j+3}=u_{j+1}u_{j+2}-u_j
\end{equation}
and $u_j^2\in\Z_{>0}$ for every $j$.  Then $u_j=2$ for every $j$.
\end{lemma}

\begin{proof}
Choose a maximal term $M=u_0$.  From
\[
 u_{-1}+u_2=u_0u_1,\qquad
 u_{-2}+u_1=u_{-1}u_0
\]
and maximality we obtain $u_{-1},u_1\le2$.  Hence each of these two neighbors belongs to
\[
 \{1,\sqrt2,\sqrt3,2\}.
\]
If both are $1$, then $u_2=M-1$ and the next term is
\[
 u_3=u_1u_2-u_0=-1,
\]
a contradiction.  Otherwise the inequalities
\[
 M(u_1-1)\le u_{-1},\qquad
 M(u_{-1}-1)\le u_1
\]
bound $M<5$.  Since each term has a positive integer square, only finitely many triples $(u_{-1},M,u_1)$ remain.  Requiring the squares of
\[
 Mu_1-u_{-1},\qquad Mu_{-1}-u_1
\]
to be integers, and requiring both terms to be positive and at most $M$, leaves the four triples listed in
\ref{app:cores}. The triples $(\sqrt2,2,\sqrt2)$ and
$(\sqrt3,2,\sqrt3)$ reach zero after two and three further terms,
respectively. The remaining nonconstant possibility is
\[
 (\ldots,2\sqrt2,\sqrt2,3,\sqrt2,2\sqrt2,1,\ldots).
\]
Two further applications of \eqref{eq:sqrt-recurrence} give $\sqrt2$ and then $-\sqrt2$, again contradicting positivity.  Thus $M=2$, and the recurrence forces every term to equal $2$.
\end{proof}

\begin{theorem}\label{thm:no-odd-period}
If $T^rq=q$ for a positive odd integer $r$, then $Tq=q$.  In particular, no positive exponent state has exact period $3$ under $T$; an exponent state returning after nine coordinate jumps already returns after three.
\end{theorem}

\begin{proof}
Proposition~\ref{prop:AB-invariants} gives $\mathcal A=\mathcal B$.  Since
\[
 \mathcal A\mathcal B=XYZ,
\]
we have $\mathcal A=\mathcal B=\sqrt{XYZ}$.  Put
\[
 U=\sqrt X,\qquad V=\sqrt Y,\qquad W=\sqrt Z.
\]
When the pair $(q_3,q_6)$ is mutated,
\[
 Z'=XY-2\sqrt{XYZ}+Z=(UV-W)^2.
\]
The positive exponent region gives $UV>W$, so the positive square root transforms by the Vieta reflection
\[
 W\longmapsto UV-W.
\]
Following the successive opposite-pair mutations therefore produces a positive periodic sequence satisfying \eqref{eq:sqrt-recurrence}, whose terms have integer squares.  Lemma~\ref{lem:sqrt-recurrence} yields
\[
 X=Y=Z=4.
\]
Now $\mathcal A=\mathcal B$ is equivalent to $q_1q_3=2q_2$, and together with the three opposite products equal to $4$ this gives
\[
 2q_i=q_{i-1}q_{i+1}
\]
for every $i$.  Hence $q$ is one of the fixed states of Theorem~\ref{thm:fixed-exponents}.
\end{proof}

\paragraph{Comparison with acyclic cluster algebras.}
For $\alpha=\beta=\gamma=1$, the exchange polynomials in
\eqref{eq:ABC} determine, up to overall sign, the candidate ordinary
cluster exchange matrix
\[
 \mathsf B(p)\defeq
 \begin{pmatrix}
 0&-c&e\\
 a&0&-f\\
 -b&d&0
 \end{pmatrix},
 \qquad p=(a,b,c,d,e,f).
\]
A matrix $\mathsf B=(b_{ij})$ is \emph{skew-symmetrizable} if
$D\mathsf B$ is skew-symmetric for some diagonal matrix $D$ with
positive integer diagonal entries. It is \emph{acyclic} if the directed
graph with an edge $i\to j$ whenever $b_{ij}>0$ has no directed cycle.
The ordinary matrix mutation at index $k$ is~\cite{FominZelevinsky2002}
\[
 (\mu_k\mathsf B)_{ij}\defeq
 \begin{cases}
  -b_{ij},& i=k\ \text{or}\ j=k,\\[2pt]
  b_{ij}+\dfrac{|b_{ik}|b_{kj}+b_{ik}|b_{kj}|}{2},& i,j\ne k.
 \end{cases}
\]
The indices $1,2,3$ correspond to jumps in $x,y,z$, respectively.

\begin{proposition}\label{prop:cluster-periodic}
For a positive integer exponent state $p=(a,b,c,d,e,f)$, the matrix
$\mathsf B(p)$ is skew-symmetrizable if and only if $ade=bcf$.
If $T^rp=p$ for some integer $r\ge1$ and every intermediate exponent state is
positive, this equality holds if and only if $Tp=p$.
For each of the seven fixed exponent states,
\[
 \mu_i\mathsf B(p)=-\mathsf B(p)
 \qquad(i=1,2,3),
\]
and no sequence of matrix mutations gives an acyclic matrix.
\end{proposition}

\begin{proof}
Since $\det\mathsf B(p)=ade-bcf$, skew-symmetrizability implies
$ade=bcf$: a skew-symmetric matrix of odd order has zero determinant.
Conversely, under this equality,
$D\defeq\operatorname{diag}(ab,bc,ae)$ makes $D\mathsf B(p)$
skew-symmetric.

Now assume $T^rp=p$ and $ade=bcf$. In the cyclic coordinates
\eqref{eq:q-order}, the latter equality is $\mathcal A=\mathcal B$.
Proposition~\ref{prop:AB-invariants} preserves this equality after
every jump. With $X=ac$, $Y=be$, and $Z=df$, an $x$-jump therefore gives
\[
 Z'=\bigl(\sqrt X\sqrt Y-\sqrt Z\bigr)^2,
 \qquad \sqrt X\sqrt Y>\sqrt Z,
\]
as in the proof of Theorem~\ref{thm:no-odd-period}. The analogous
formulas for the other jumps produce a positive periodic sequence
satisfying \eqref{eq:sqrt-recurrence}, whose terms have integer squares.
Lemma~\ref{lem:sqrt-recurrence}
gives $X=Y=Z=4$, and the same argument yields
\eqref{eq:fixed-local}. Hence $Tp=p$. Conversely, if $Tp=p$,
Proposition~\ref{prop:AB-invariants} gives
$\mathcal A-\mathcal B=-(\mathcal A-\mathcal B)$, so $ade=bcf$.

For a fixed state, \eqref{eq:fixed-local} gives $ae=2f$ and $bc=2d$.
Thus mutation at $1$ replaces the entries $-f,d$ by
$ae-f=f$ and $d-bc=-d$, while reversing the signs in row and column
$1$. It follows that $\mu_1\mathsf B(p)=-\mathsf B(p)$; the other
two cases follow cyclically. Since
$\mu_i(-\mathsf B)=-\mu_i(\mathsf B)$, every matrix obtained by
mutations is either $\mathsf B(p)$ or $-\mathsf B(p)$. Both have
a directed triangle, so neither is acyclic.
\end{proof}

In particular, none of the twenty-two genuine exponent two-cycles can
be represented by the displayed exchanges in an ordinary
skew-symmetrizable rank-three cluster seed. The seven fixed exponent
states admit such seeds, but their matrices cannot be mutated to
acyclic ones. Thus Zhang's acyclic classification
\cite[Theorem~B]{Zhang2025} does not apply directly to these
three-variable exchanges.

\section{Sign changes and returns}\label{sec:signs}

For each exponent two-cycle, choose the representative $p$ listed in
the first exponent column of \ref{app:cycles}. In Sections
\ref{sec:signs}--\ref{sec:positive-lifts}, sign assignments and
realizations are counted at these chosen representatives. A sign
assignment together with its representative specifies a signed
polynomial state.

Encode the signs by elements of $\Ftwo$:
\[
 \alpha=(-1)^A,\qquad
 \beta=(-1)^B,\qquad
 \gamma=(-1)^G,
 \qquad
 (A,B,G)\in\Ftwo^3.
\]
For an exponent state $p=(a,b,c,d,e,f)$, put $p^*\defeq T(p)$ and let $b^*=ad_1-b$ be its second exponent.  Stars on $A,B,G$ denote the sign exponents after one block.

\begin{proposition}\label{prop:sign-matrix}
The sign update \eqref{eq:sign-block} is the linear map
\begin{equation}\label{eq:sign-matrix}
 \begin{pmatrix}A^*\\B^*\\G^*\end{pmatrix}
 =
 M(p)
 \begin{pmatrix}A\\B\\G\end{pmatrix},
 \qquad
 M(p)=
 \begin{pmatrix}
 1+ac & a & b^*\\
 c & 1 & 0\\
 0 & 0 & 1
 \end{pmatrix}
 \pmod2.
\end{equation}
\end{proposition}

\begin{proof}
Equation \eqref{eq:sign-block} becomes
\[
 B^*=B+cA,\qquad
 G^*=G,\qquad
 A^*=A+a(B+cA)+b^*G
\]
over $\Ftwo$.
\end{proof}

\begin{theorem}\label{thm:84-signs}
Let $p\leftrightarrow p^*$ be one of the twenty-two exponent cycles of exact period two, and let $I$ denote the $3\times3$ identity matrix.  The sign triple $(\alpha,\beta,\gamma)$ returns after the corresponding six jumps if and only if its exponent vector $(A,B,G)^{\mathsf T}$ lies in
\[
 \ker\bigl(M(p^*)M(p)-I\bigr).
\]
For these twenty-two chosen representatives, there are exactly $84$ returning sign assignments:
\begin{itemize}[leftmargin=2em]
\item four exponent cycles admit all eight sign triples;
\item eight exponent cycles admit four sign triples;
\item ten exponent cycles admit two sign triples.
\end{itemize}
In particular, every exponent two-cycle has at least one returning sign assignment.
\end{theorem}

\begin{proof}
The kernel condition follows from Proposition~\ref{prop:sign-matrix}.  Evaluating the $3\times3$ matrices for the twenty-two exponent cycles of Corollary~\ref{cor:22-cycles} gives kernels with $8$, $4$, and $2$ elements, with the stated frequencies.  \ref{app:cycles} lists every cycle and every returning sign pattern.
\end{proof}

For four cycles, $M(p^*)M(p)=I$, so all eight sign triples return after six jumps.  In the exponent order $(a,b,c,d,e,f)$, one example is
\[
 (2,2,3,2,2,3)\longleftrightarrow(4,6,1,4,1,1),
\]
the cycle used in the arithmetic example below.

\section{Integral realizations and a single equation}\label{sec:arithmetic-lifts}

The exponent and sign classifications concern polynomial states. We now ask which of them are realized by sequences of exact integer quotient jumps. An \emph{integral six-jump realization} of one of these returning signed polynomial states consists of two successive $xyz$ blocks starting from that state, with exact integer quotients and nonzero integer coordinates throughout. It is \emph{positive} if the initial and all subsequent coordinates are positive. The polynomial state returns after the six jumps, while the initial and final coordinate triples may differ.

\begin{proposition}\label{prop:all84-integral}
Every one of the $84$ returning signed polynomial states in Theorem~\ref{thm:84-signs} has a pairwise-coprime nonzero integer six-jump realization.
\end{proposition}

\begin{proof}
The seventy-four positive cases of Theorem~\ref{thm:positive-classification} below already give such realizations.  For each of the remaining ten sign states, direct substitution shows that one of the four pairwise-coprime starting triples
\[
 (-1,-1,-1),\qquad
 (-1,-1,1),\qquad
 (-1,1,-1),\qquad
 (-1,-2,-1)
\]
produces six nonzero integral quotient jumps.  Lemma~\ref{lem:coprime-preserved} keeps the orbit pairwise coprime.
\end{proof}

For a fixed signed polynomial state \eqref{eq:ABC}, define
\begin{equation}\label{eq:FABC}
 F_{A,B,C}(x,y,z)
 \defeq
 yz\,A(y,z)+xz\,B(x,z)+xy\,C(x,y).
\end{equation}
The term \emph{binomial} refers to $A,B,C$ separately; the combined
polynomial $F_{A,B,C}$ has six distinct monomial terms.

\begin{proposition}\label{prop:F-equivalence}
At an integer point with nonzero, pairwise-coprime coordinates,
\[
 x\mid A(y,z),\qquad
 y\mid B(x,z),\qquad
 z\mid C(x,y)
\]
if and only if
\begin{equation}\label{eq:F=kxyz}
 F_{A,B,C}(x,y,z)=kxyz
\end{equation}
for some $k\in\Z$.  In that case
\begin{equation}\label{eq:k-sum}
 k=
 \frac{A(y,z)}x+
 \frac{B(x,z)}y+
 \frac{C(x,y)}z.
\end{equation}
\end{proposition}

\begin{proof}
The forward implication follows from \eqref{eq:k-sum}.  Conversely, reduce \eqref{eq:F=kxyz} modulo $x$.  Pairwise coprimality gives $\gcd(x,yz)=1$, hence $x\mid A(y,z)$.  Reducing modulo $y$ and $z$ proves the other two divisibilities.
\end{proof}

\begin{corollary}\label{cor:allplus}
Every one of the twenty-two exponent two-cycles has a positive pairwise-coprime arithmetic realization: choose the sign triple $(\alpha,\beta,\gamma)=(1,1,1)$ and the initial triple
\[
 (x,y,z)=(1,1,1).
\]
\end{corollary}

\begin{proof}
The sign triple $(1,1,1)$ always returns.  Every numerator in every quotient jump is then a sum of positive powers, so positivity is automatic.  The three initial divisibilities hold because every coordinate is $1$.  Reciprocal transport and Lemma~\ref{lem:coprime-preserved} preserve integrality at each jump.
\end{proof}

The next example returns to the same polynomial $F_{A,B,C}$ with a different integer $k$ in \eqref{eq:F=kxyz}.
\begin{example}\label{ex:k-change}
Take
\[
 A=y^2+z^2,\qquad
 B=x^3-z^2,\qquad
 C=x^2+y^3.
\]
The exponent state
\[
 (2,2,3,2,2,3)
\]
belongs to the two-cycle
\[
 (2,2,3,2,2,3)
 \longleftrightarrow
 (4,6,1,4,1,1),
\]
and the sign triple $(+,-,+)$ returns after six jumps.  Starting at
\[
 (x,y,z)=(26,11,3)
\]
the six quotient jumps end at
\[
 (197,47,47544),
\]
while the polynomial state returns exactly.  For the single polynomial
\[
 F=yz(y^2+z^2)+xz(x^3-z^2)+xy(x^2+y^3),
\]
the initial and terminal points satisfy
\[
 F=2271\,xyz
\]
and
\[
 F=-36457341\,xyz,
\]
respectively.  The initial and final triples therefore satisfy equations with the same polynomial $F$ but different integer values of $k$.
\end{example}

\section{Positive integral realizations}\label{sec:positive-lifts}

We now determine which of the $84$ returning signed polynomial states admit positive pairwise-coprime integral six-jump realizations.
The constructions below also illustrate the multivariable use of the
Mordell--Schinzel method of constructing solutions backwards.  Once a polynomial-state
sequence is fixed, one may prescribe terminal coordinates and solve its
quotient equations in reverse order.  Lemma~\ref{lem:block-predecessor}
gives the three-variable predecessor formulas, and the constructions
$\mathsf P(u)$ and $\mathsf R(X,Y,z)$ below choose suitable terminal or
intermediate data.  Exact divisibility, positivity, and the required
coprimality conditions are part of these constructions, not consequences
of choosing a terminal tuple alone.  Their purpose here is to establish
finite six-jump realizations; indefinite repetition is a separate question.

For one block, write the three successive quotient coordinates as $(X,Y,Z)$.  With the notation of Section~\ref{sec:binomial}, the forward equations are
\begin{equation}\label{eq:block-coordinate}
 xX=y^a+\alpha z^b,\qquad
 yY=X^c+\beta_1z^{d_1},\qquad
 zZ=X^{e_2}+\gamma Y^{f_1},
\end{equation}
where $\beta_1=\beta\alpha^c$.

\Needspace{6\baselineskip}
We first recover the predecessor of one complete $xyz$ block.
\begin{lemma}\label{lem:block-predecessor}
Fix a signed polynomial state in the positive exponent region.  Given a positive terminal triple $(X,Y,Z)$, a positive predecessor is uniquely determined, when it exists, by
\begin{equation}\label{eq:predecessor}
 z=\frac{X^{e_2}+\gamma Y^{f_1}}{Z},
 \qquad
 y=\frac{X^c+\beta_1z^{d_1}}{Y},
 \qquad
 x=\frac{y^a+\alpha z^b}{X}.
\end{equation}
Thus, if the three displayed fractions are positive integers, these formulas give an explicit positive integral predecessor for the block.
\end{lemma}

\begin{proof}
This is \eqref{eq:block-coordinate} solved in reverse order.
\end{proof}

For one construction, require the triple after the first block to be $(X,Y,1)$.  Stars refer to the polynomial state reached after those three jumps; intermediate starred quantities are defined by the same formulas as in Section~\ref{sec:binomial}.  Choose a positive integer $u$, require the first two coordinates after the second block to be $u$ and $1$, and denote the third by $W$.  Its three quotient equations then require
\begin{equation}\label{eq:P-template}
 Y=u^{c^*}+\beta_1^*,\qquad
 X=\frac{Y^{a^*}+\alpha^*}{u},\qquad
 W=u^{e_2^*}+\gamma^*.
\end{equation}
Whenever these are positive integers and Lemma~\ref{lem:block-predecessor} recovers a positive integral predecessor for the first block, we call the resulting construction $\mathsf P(u)$.

For a second construction, denoted $\mathsf R(X,Y,z)$, prescribe positive integers $X,Y,z$, where $X,Y$ are the coordinates after the first two jumps and $z$ is the initial third coordinate.  Set
\[
 Z=\frac{X^{e_2}+\gamma Y^{f_1}}z
\]
and use \eqref{eq:predecessor} to recover $y$ and $x$.

\begin{lemma}\label{lem:midpoint-congruence}
Consider three successive positive integer quotient jumps starting from a positive integer triple with pairwise-coprime coordinates, and assume $\beta_1=-1$ and $ad_1\ge b$.  Then
\begin{equation}\label{eq:midpoint-congruence}
 X\mid (-1)^a z^{ad_1-b}+\alpha Y^a.
\end{equation}
Moreover positivity of the second jump gives
\begin{equation}\label{eq:second-ineq}
 X^c>z^{d_1}.
\end{equation}
If $\gamma=-1$, positivity of the third jump also gives
\begin{equation}\label{eq:third-ineq}
 X^{e_2}>Y^{f_1}.
\end{equation}
\end{lemma}

\begin{proof}
Modulo $X$, the first equation of \eqref{eq:block-coordinate} gives
\[
 y^a\equiv-\alpha z^b.
\]
The second gives
\[
 yY\equiv-z^{d_1}.
\]
Raise the latter congruence to the $a$th power and substitute the former.  Pairwise coprimality permits cancellation of $z^b$, yielding \eqref{eq:midpoint-congruence}.  The numerator of the second jump is positive, as is that of the third jump when $\gamma=-1$; rearranging these inequalities gives \eqref{eq:second-ineq} and \eqref{eq:third-ineq}.
\end{proof}

\begin{theorem}\label{thm:positive-classification}
Of the $84$ returning signed polynomial states, exactly $74$ admit a positive pairwise-coprime integral six-jump realization.  The remaining ten are the sign patterns in the last column of \ref{app:cycles}.
\end{theorem}

\begin{proof}
For existence, \ref{app:positive-certificates} assigns one of the six constructions
\[
 \mathsf P(2),\quad \mathsf P(3),\quad
 \mathsf R(2,1,1),\quad
 \mathsf R(4,3,1),\quad
 \mathsf R(5,2,1),\quad
 \mathsf R(13,1,6)
\]
to $67$ signed states.  The initial triples for the seven remaining states are given in Table~\ref{tab:residual-seven}.  In every case the three divisions in \eqref{eq:predecessor} and the three divisions in the second block are exact and positive.

For nonexistence, suppose that one of the ten excluded states in
\ref{app:cycles} has a positive pairwise-coprime six-jump realization.
Write its points after three and six jumps as $(X,Y,Z)$ and $(u,v,w)$,
respectively.  In particular, $\gcd(Y,z)=1$, by
Lemma~\ref{lem:coprime-preserved}.  Every excluded state has
$\beta_1=-1$, so Lemma~\ref{lem:midpoint-congruence} applies.
The four cases below cover all ten exclusions in \ref{app:cycles}.

First consider the five entries in the following table.  A row number
refers to the exponent cycle in \ref{app:cycles}.  The column
headed $G$ gives an integer divisible by $X$; the last two columns are
strict inequalities forced by positivity of the second, third, and
fourth jumps.
\[
\begin{array}{c|c|c|c|c}
 \text{row}&(\alpha,\beta,\gamma)&G&\text{bound on }X&
 \text{comparison}\\ \hline
 5 &(+,-,+)&Y-z       &X>z  &Y<z\\
 11 &(+,-,+)&Y-z^2    &X>z^2&Y<z^2\\
 16 &(-,+,+)&z-Y^2    &X>z  &Y^2<z\\
 18 &(-,+,+)&z-Y^2    &X>z  &Y^2<z\\
 19 &(-,+,+)&z^2-Y^2  &X>z^2&Y<z
\end{array}
\]
Here is the verification of the inequalities involving $Y$ and $z$.
For row~5, the third and fourth jumps give
\[
 zZ=X^2+Y^4>Y^4,\qquad Y^3>Z,
\]
so $z>Y$.  For row~11 they give
\[
 zZ=X+Y^2>Y^2,\qquad Y^3>Z^2,
\]
and hence $Y^3>Y^4/z^2$, so $z^2>Y$.
For rows~16 and~18 they give, respectively,
\[
 zZ=X^2+Y^6>Y^6,\quad Y^4>Z,
 \qquad
 zZ=X^2+Y^5>Y^5,\quad Y^3>Z,
\]
so $z>Y^2$ in both cases.  Finally, row~19 gives
\[
 zZ=X+Y^3>Y^3,\qquad Y^4>Z^2,
\]
so $Z<Y^2$ and therefore $z>Y$.
The bounds on $X$ follow from the second-jump inequalities
$X^2>z^3$, $X^2>z^6$, $X>z^2$, $X>z^3$, and $X>z^4$,
respectively.  Thus every row gives $0<|G|<X$, contradicting $X\mid G$.

Second, consider row~11 with signs $(+,-,-)$ and row~19 with signs
$(-,+,-)$.  Their divisibility conditions from Lemma~\ref{lem:midpoint-congruence} and positivity bounds are
\[
\begin{array}{c|c|c}
 \text{row and signs}&\text{divisibility condition}&\text{bounds}\\ \hline
 11, (+,-,-)&X\mid Y-z^2&X>z^3,\quad X>Y^2\\
 19, (-,+,-)&X\mid z^2-Y^2&X>z^4,\quad X>Y^3.
\end{array}
\]
Both terms in each difference are strictly less than $X$; divisibility
therefore forces $Y=z^2$ in the first case and $Y=z$ in the second.
Since $\gcd(Y,z)=1$, both cases have $Y=z=1$.  The third jump then gives
$Z=X-1\ge1$, whereas the numerator of the fourth jump is $1-Z^2$.
It cannot be positive.

Third, consider rows~17 and~19, both with signs $(+,-,-)$.  In both
cases the divisibility condition is
\[
 X\mid z^2+Y^2.
\]
Positivity gives $X>z^3$ and $X>Y^4$ in row~17, and $X>z^4$ and
$X>Y^3$ in row~19.  If $M\defeq\max(z,Y)\ge2$, then
\[
 X>M^3\ge2M^2\ge z^2+Y^2,
\]
which contradicts divisibility.  Hence $z=Y=1$, and the bounds and
divisibility force $X=2$.  The third jump gives $Z=1$; the fourth gives
$u=1$.  The fifth-jump numerator is $u^2-Z$ in row~17 and $u^3-Z^2$
in row~19, so it is zero in both cases.

Finally, row~11 with signs $(-,-,-)$ gives
\[
 X\mid z^2+Y,\qquad X>z^3,\qquad X>Y^2.
\]
For positive integers $(z,Y)\ne(1,1)$,
\[
 \max(z^3,Y^2)\ge z^2+Y.
\]
Indeed, if $Y\le z$ and $z\ge2$, then $z^3\ge z^2+z\ge z^2+Y$;
if $Y\ge z+1$, then $Y^2-Y\ge z(z+1)>z^2$.
Thus $X>z^2+Y$, again impossible.  The remaining case has $z=Y=1$
and $X=2$.  It gives $Z=u=1$, and the fifth-jump numerator
$u^4-Z^2$ vanishes.  This excludes the last excluded state and
completes the classification.
\end{proof}

\begin{remark}\label{rem:one-block-positive}
Theorem~\ref{thm:positive-classification} concerns existence of one positive six-jump return of the polynomial state.  It does not assert that the coordinate orbit remains positive under arbitrarily many repetitions of that six-jump block.
\end{remark}

\section{Reverse jumps, roots, and infinite descent}\label{sec:reverse-descent}

We use decreasing height to reach a root or to prove nonexistence by infinite descent. The vertices must record both their coordinates and all polynomial data needed to determine the permitted jumps. When a jump is reversible on this full state, reversing a decreasing jump gives an increasing one.

\subsection{A general descent principle}\label{subsec:general-root-descent}

Fix the number $r$ of coordinates and a specified jump procedure.  Its vertices are
\[
 \mathcal V\defeq
 \{(\sigma,\bx):\bx\in(\Z\setminus\{0\})^r
       \text{ is an admissible solution in the state }\sigma\}.
\]
The state $\sigma$ records all information needed to determine the permitted jumps: the polynomials, any weights or auxiliary quotients, the normalization rule, and the current position in a prescribed sequence of coordinate jumps when one is used. An edge records a permitted jump and its resulting state. Thus the same coordinate tuple can occur at different vertices. We allow directed edges. A reverse jump must restore the full vertex, including $\sigma$. Throughout this subsection, a component means a connected component after edge directions are ignored; a path must follow the permitted directions.

For $v=(\sigma,\bx)$, let $a_1\le\cdots\le a_r$ be the absolute coordinate values in nondecreasing order, and define
\[
 H(v)\defeq(a_r,\ldots,a_1),
 \qquad
 h(v)\defeq\sum_{i=1}^r|x_i|.
\]
We compare $H$ lexicographically: at the first position where two height vectors differ, the smaller entry gives the smaller height. Both heights depend only on the coordinates, while admissibility depends on $\sigma$. On a fixed positive solution set $\mathcal S$ preserved by involutions, this gives the usual solution graph; when $\sigma$ is fixed we also write $H(\bx)$.

\begin{definition}\label{def:procedure-root}
A vertex is a \emph{local root}, or simply a \emph{root}, if no permitted jump decreases $H$.  A jump is \emph{neutral} if it preserves $H$.  A root is \emph{stable under neutral transitions} if every vertex reachable from it by a finite sequence of neutral jumps is a local root.  A \emph{component minimum} is a vertex whose height is least in its component.
\end{definition}

We call a path \emph{nonincreasing} if $H$ never increases along it. A root is stable under neutral transitions exactly when no such path starting at that root reaches a smaller height: all jumps before the first strict decrease would be neutral. In particular, the definitions give
\[
 \text{component minimum}
 \ \Longrightarrow\
 \text{root stable under neutral transitions}
 \ \Longrightarrow\
 \text{local root}.
\]
The examples below show that neither implication can be reversed in general.

\begin{proposition}\label{prop:root-descent}
Every vertex admits a finite sequence of strict descents ending at a local root.  If each jump changes at most one coordinate, every such sequence starting at $v_0=(\sigma_0,\bx^{(0)})$ has length at most
\[
 h(v_0)-r=\sum_{i=1}^r|x_i^{(0)}|-r.
\]
Every vertex also admits a finite nonincreasing path to a root stable under neutral transitions.  Every nonempty component contains a component minimum, and hence a root.  In particular, if the root set is empty, then $\mathcal V$ is empty.
\end{proposition}

\begin{proof}
Repeatedly choose a height-decreasing jump.  The lexicographic order on $\Z_{>0}^r$ is well founded, so the process terminates at a root.  If no root exists, the assumption that a starting solution exists would instead produce an infinite strictly decreasing sequence of positive integer height vectors.

Suppose an edge $v\to v'$ replaces only $x_i$ by $x_i'$.  Decreasing one absolute coordinate cannot increase any entry after the list is sorted in decreasing order, and at least one entry strictly decreases. Increasing that coordinate has the opposite effect.  Consequently
\[
 H(v')<H(v)
 \quad\Longleftrightarrow\quad
 |x_i'|<|x_i|
 \quad\Longleftrightarrow\quad
 h(v')<h(v),
\]
and equality of $H$ is equivalent to equality of $h$ on this edge.  Each strict descent therefore decreases the integer $h$ by at least one.  Since every coordinate is nonzero, $h\ge r$, proving the bound.

Let $\mathcal R_{\le}(v)$ be the set of vertices reachable from $v$ by finite nonincreasing paths, including the path of length zero.  Choose $w\in\mathcal R_{\le}(v)$ of least height.  A neutral path from $w$ followed by a strict descent would give a vertex of smaller height in $\mathcal R_{\le}(v)$, a contradiction.  Hence $w$ is stable under neutral transitions.  Finally, the same well-ordering of height vectors gives a least height in every nonempty component.
\end{proof}

For an admissible quotient jump
\[
 x_i'=\frac{Q_{i,\sigma}(\widehatx{i})}{x_i},
\]
where $Q_{i,\sigma}$ denotes the chosen numerator, the strict-descent test is
\[
 H(v')<H(v)
 \quad\Longleftrightarrow\quad
 0<|Q_{i,\sigma}(\widehatx{i})|<x_i^2.
\]
The new coordinate must be an integer and nonzero, and the updated polynomial relations must hold; these conditions are part of admissibility. The inequality alone does not establish them. The equivalence with $h$ used in the proof applies to a single coordinate jump, not to arbitrary pairs of vertices.

The states $\sigma$ need not form a finite set or repeat periodically.  If the jumps are reversible, reversing a strict descent constructs the original vertex from a local root by strict ascents.  Every admissible vertex with coordinates $(\pm1,\ldots,\pm1)$ is a component minimum.  In a directed graph, however, a component minimum need not be reachable from a specified starting vertex; the nonincreasing path in Proposition~\ref{prop:root-descent} has the weaker conclusion stated there.

The termination bound controls strict descents, not the number of neutral transitions.  Even at fixed coordinates there may be infinitely many possible states $\sigma$.  Thus existence of a root stable under neutral transitions does not by itself give a finite procedure for recognizing one.  Nor does the root principle imply finiteness or uniqueness of roots, or return to a prescribed equation.  For the fixed generalized Hurwitz equations treated in~\cite[Definition~3.3 and Theorem~3.5]{FineEtAl2016Hurwitz}, the definition by absence of a lower-height solution obtainable by one Vieta jump and the uniqueness theorem for each orbit are separate steps.

\paragraph{A neutral jump can permit a later decrease.}
The distinction already occurs in two variables. For integer polynomials $A,B,P$ with $P\ne0$ and $c\in\Z$, write
\[
 \mathcal S(A,B):\qquad y\mid A(x),\qquad x\mid B(y),
 \qquad
 \mathcal R_cP(T)\defeq T^{\deg P}P(c/T).
\]
We denote a realizing vertex by $(A,B;x,y)$.  In each example the only permitted moves are the two coordinate quotient jumps, with the state transport specified there.  For transport without normalization, a $y$-jump replaces $B$ by $\mathcal R_{A(0)}B$, and an $x$-jump replaces $A$ by $\mathcal R_{B(0)}A$, as in Theorem~\ref{thm:face-transport}.

Take
\[
 A(T)\defeq T+1,\qquad B(T)\defeq4T+1,\qquad C(T)\defeq T+4.
\]
At $(A,B;3,2)$ the coordinate quotients are $A(3)/2=2$ and $B(2)/3=3$, so this vertex is a local root.  Nevertheless,
\[
 \mathcal R_1B=C,\qquad \mathcal R_4A=C
\]
give the admissible path
\[
 (A,B;3,2)\longrightarrow(A,C;3,2)\longrightarrow(C,C;2,2),
\]
with heights
\[
 (3,2)\longrightarrow(3,2)\longrightarrow(2,2).
\]
The first jump is neutral but changes the polynomial state; the next quotient is $C(2)/3=2$ and decreases the height. Thus the initial root is not stable under neutral transitions. The final coordinates need not be coprime, since Theorem~\ref{thm:face-transport} has no such hypothesis. The reverse $x$-quotient returns the coordinates to $(3,2)$ but gives $\mathcal R_4C=4A$, so it does not restore the preceding full vertex.

\paragraph{Different descent branches can reach different roots.}
For this example choose the following normalization procedure. After each coordinate quotient jump, divide the transformed polynomial by its content and choose its sign so that its leading coefficient is positive. We permit this step only if the transformed polynomial is nonzero and its content is coprime to the unchanged coordinate; division by the content then preserves the required divisibility. The weighted construction of Theorem~\ref{thm:weighted-primitive-transport} retains this content as a weight instead.

Put
\[
 \begin{aligned}
 A(T)&\defeq T^2-7,&
 B(T)&\defeq5T^2+17T-7,\\
 C(T)&\defeq T^2+17T-35,&
 D(T)&\defeq T^2-175.
 \end{aligned}
\]
The reciprocal identities
\[
 \begin{aligned}
 \mathcal R_{-7}B&=-7C,&
 \mathcal R_{-35}A&=-7D,\\
 \mathcal R_{-7}C&=-7B,&
 \mathcal R_{-35}D&=-175A
 \end{aligned}
\]
give a path with admissible reverse edges:
\[
 (A,B;3,1)\longleftrightarrow(A,C;3,2)
              \longleftrightarrow(D,C;1,2).
\]
For the two forward edges the quotients are $A(3)/1=2$ and $C(2)/3=1$, and the discarded contents are both $7$, coprime to the unchanged coordinates $3$ and $2$.  In reverse order the discarded contents are $175$ and $7$, again coprime to the unchanged coordinates $2$ and $3$.

Both endpoints are local roots.  At the first, the possible replacement coordinates are
\[
 x'=\frac{B(1)}3=5,\qquad y'=\frac{A(3)}1=2;
\]
at the last, they are
\[
 x'=\frac{C(2)}1=3,\qquad y'=\frac{D(1)}2=-87.
\]
Every move increases the height. The two outgoing moves not already checked on the displayed path give the following polynomials before normalization:
\[
 \mathcal R_{-7}A=-7A,\qquad
 \mathcal R_{-175}C=-35(T^2+85T-875);
\]
their unchanged coordinate is $1$, so these cancellations are admissible as well.  If only positive coordinates are allowed, the last $y$-move is excluded and the same root conclusions hold.  Neither endpoint admits a neutral first jump, so both roots are stable under neutral transitions.  Their heights along the path are
\[
 H:\ (3,1)\longrightarrow(3,2)\longrightarrow(2,1),
 \qquad
 h:\ 4\longrightarrow5\longrightarrow3.
\]
The first endpoint is therefore not a component minimum.  From the middle vertex, the two choices of strict descent reach different roots.

Along this path, the corresponding equations have the same parameter $k=8$:
\begin{center}
\begin{tabular}{@{}lc@{}}
\toprule
Polynomial defining $F_j(x,y)=8xy$ & Realizing point\\
\midrule
$F_0(x,y)\defeq x^2+5y^2+17y-7$ & $(3,1)$\\
$F_1(x,y)\defeq5x^2+y^2+17y-35$ & $(3,2)$\\
$F_2(x,y)\defeq x^2+5y^2+85y-175$ & $(1,2)$\\
\bottomrule
\end{tabular}
\end{center}
The preservation is explained by the exact identities
\[
 \begin{aligned}
 yF_1(x,z)-zF_0(x,y)&=(5y-z)(A(x)-yz),\\
 xF_2(u,y)-uF_1(x,y)&=(5x-u)(C(y)-xu).
 \end{aligned}
\]
Thus the relations defining the jumps, $yz=A(x)$ and $xu=C(y)$, transport the respective equations without changing $k$.

These examples show that decreasing height on the graph of states and solutions differs from finding a minimum among solutions of a fixed equation. We now consider a four-variable equation preserved by every coordinate jump, for which explicit bounds give finite classifications of roots.

\subsection{A four-variable family with cross terms}

For $\lambda\in\Z_{>0}$ define
\begin{equation}\label{eq:Qlambda}
 Q_\lambda(x,y,z,w)
 \defeq
 x^2+y^2+z^2+w^2
 +\lambda(xy+xz+xw+yz+yw+zw).
\end{equation}
We seek positive integer solutions of
\begin{equation}\label{eq:Qlambda-level}
 Q_\lambda(x,y,z,w)=kxyzw,
 \qquad k\in\Z_{>0}.
\end{equation}

As a quadratic in $w$, equation~\eqref{eq:Qlambda-level} is
\begin{equation}\label{eq:Qlambda-quadratic}
 w^2-
 \bigl(kxyz-\lambda(x+y+z)\bigr)w
 +Q_\lambda(x,y,z)=0,
\end{equation}
where $Q_\lambda(x,y,z)$ denotes the same quadratic expression restricted to three variables.

\begin{proposition}\label{prop:Qlambda-jump}
For every positive integral solution of~\eqref{eq:Qlambda-level}, the coordinate jump
\begin{equation}\label{eq:w-jump}
 w'
 \defeq
 kxyz-\lambda(x+y+z)-w
 =\frac{Q_\lambda(x,y,z)}{w}
\end{equation}
gives another positive integral solution when $w$ is replaced by $w'$, and applying the same jump again restores $w$. The analogous statement holds in each coordinate. Moreover
\begin{equation}\label{eq:Phi-lambda}
 \Phi_\lambda(x,y,z,w)
 \defeq
 \frac{Q_\lambda(x,y,z,w)}{xyzw}
\end{equation}
is unchanged by each coordinate jump whose initial and final coordinates are all nonzero rational numbers.
\end{proposition}

\begin{proof}
Equation~\eqref{eq:Qlambda-quadratic} shows that $w$ and $w'$ are the two roots of the quadratic; the product relation gives $ww'=Q_\lambda(x,y,z)>0$, so $w'$ is a positive integer and the transformation is involutive.  For invariance of~\eqref{eq:Phi-lambda}, put $q=Q_\lambda(x,y,z)$. On the stated domain, $w'=q/w\ne0$, and
\[
 Q_\lambda(x,y,z,w')
 =\frac{q}{w^2}Q_\lambda(x,y,z,w),
\]
while $xyzw'=xyzq/w$, so the quotient~\eqref{eq:Phi-lambda} is unchanged.
\end{proof}

Preservation under coordinate jumps connects~\eqref{eq:Qlambda-level} with Diophantine equations studied through cluster transformations~\cite{Lampe2016,BaoLi2025} and higher-dimensional Vieta mutation graphs~\cite{Cho2026}. The usual generalized Hurwitz equation has no pairwise cross terms~\cite{FineEtAl2016Hurwitz}, whereas the family~\eqref{eq:Qlambda} contains every such term.

\subsection{Bounds for roots}

Sort a positive solution as
\[
 1\le x\le y\le z\le w.
\]
Call it a root if the replacement~\eqref{eq:w-jump} of a largest coordinate satisfies $w'\ge w$.

\begin{lemma}\label{lem:largest-descent}
If $w'<w$, replacing one largest coordinate $w$ by $w'$ and resorting strictly decreases $H$.  If $w'\ge w$, then no coordinate jump decreases $H$.
\end{lemma}

\begin{proof}
The first assertion follows from the definition of $H$. For a smaller coordinate, say $x\le w$, the replacement is
\[
 x'=\frac{Q_\lambda(y,z,w)}x.
\]
Since $Q_\lambda(y,z,w)>w^2$, one has $x'>w^2/x\ge w$.  The same argument applies to $y$ and $z$.  Thus at a root no smaller-coordinate jump can decrease the height either.
\end{proof}

\begin{theorem}\label{thm:root-strip}
Let $1\le x\le y\le z\le w$ be a root of~\eqref{eq:Qlambda-level} and put
\[
 n\defeq kxy.
\]
Then
\begin{equation}\label{eq:root-strip}
 \lambda+3\le n\le6\lambda+4.
\end{equation}
More precisely, with
\[
 A\defeq n-\lambda-2,
 \qquad
 C\defeq x^2+y^2+\lambda xy,
\]
one has
\begin{equation}\label{eq:root-lower-z}
 Az\ge\lambda(x+y)
\end{equation}
and
\begin{equation}\label{eq:root-upper-z}
 Az^2\le C+2\lambda(x+y)z.
\end{equation}
Consequently, for every fixed $\lambda$ there are only finitely many roots and only finitely many positive integers $k$ for which~\eqref{eq:Qlambda-level} can have a positive integral solution.
\end{theorem}

\begin{proof}
Let
\[
 f(t)=t^2-\bigl((n-\lambda)z-\lambda(x+y)\bigr)t+Q_\lambda(x,y,z),
\]
so that its two roots are $w\le w'$.  Since $z\le w\le w'$,
\[
 w+w'=(n-\lambda)z-\lambda(x+y)\ge2z,
\]
which is~\eqref{eq:root-lower-z}; in particular $A>0$, hence $n\ge\lambda+3$.

Also $f(z)=(z-w)(z-w')\ge0$.  Expanding $f(z)$ gives exactly
\[
 C+2\lambda(x+y)z-Az^2\ge0,
\]
which is~\eqref{eq:root-upper-z}.  Divide by $z^2$ and use $x/z,y/z\le1$:
\[
 A
 \le
 \frac{x^2+y^2+\lambda xy}{z^2}
 +2\lambda\frac{x+y}{z}
 \le 2+5\lambda.
\]
Thus $n=A+\lambda+2\le6\lambda+4$.

Finally $xy\mid n$, so the bound on $n$ leaves finitely many pairs $(x,y)$.  For each pair,~\eqref{eq:root-upper-z} bounds $z$, and then~\eqref{eq:Qlambda-quadratic} leaves at most two values of $w$.
\end{proof}

The proof gives an exact enumeration procedure.  For each $n$ in~\eqref{eq:root-strip}, enumerate $x\le y$ with $xy\mid n$, set $k=n/(xy)$, and enumerate the finitely many values of $z\ge y$ satisfying~\eqref{eq:root-lower-z}--\eqref{eq:root-upper-z}.  Put
\[
 S=(n-\lambda)z-\lambda(x+y),
 \qquad
 \Delta=S^2-4Q_\lambda(x,y,z).
\]
A root occurs precisely when $\Delta$ is a square of the same parity as $S$ and
\[
 w=\frac{S-\sqrt\Delta}{2}\ge z.
\]
This is the exact arithmetic check used below and in \ref{app:Qlambda-roots}.

\subsection{The case \texorpdfstring{$\lambda=1$}{lambda=1}}

For $\lambda=1$,
\[
 Q_1(x,y,z,w)=h_2(x,y,z,w),
\]
the complete homogeneous symmetric polynomial of degree two.  Equivalently,
\[
 2Q_1=(x+y+z+w)^2+x^2+y^2+z^2+w^2.
\]
Let $e_1,\ldots,e_5$ be the standard coordinate vectors in $\mathbb R^5$. The $A_4$ root lattice consists of the integer vectors whose coordinates sum to zero. The displayed identity identifies $2Q_1$ with its squared Euclidean norm in the basis $e_i-e_5$, $1\le i\le4$.

\begin{theorem}\label{thm:lambda1-classification}
The equation
\begin{equation}\label{eq:lambda1}
 x^2+y^2+z^2+w^2+xy+xz+xw+yz+yw+zw=kxyzw
\end{equation}
has a positive integral solution if and only if
\[
 k\in\{2,8,10\}.
\]
Every positive solution descends to one of the four sorted roots
\[
 \begin{array}{c|c}
 k&\text{roots}\\
 \hline
 2&(1,2,7,7),\ (2,2,2,4)\\
 8&(1,1,1,2)\\
 10&(1,1,1,1).
 \end{array}
\]
For every other positive integer $k$, nonexistence follows by infinite descent.
\end{theorem}

\begin{proof}
Theorem~\ref{thm:root-strip} gives $4\le n=kxy\le10$. Enumerating $x\le y$ with $xy\mid n$ and $z\ge y$ satisfying~\eqref{eq:root-lower-z}--\eqref{eq:root-upper-z} produces $38$ candidate quadruples $(k,x,y,z)$. The exact discriminant test above leaves precisely the four roots displayed. Substitution verifies each displayed root, so the three listed values of $k$ occur.

For any other positive solution, Lemma~\ref{lem:largest-descent} gives a strictly smaller sorted height until a root is reached.  Hence a solution for $k\notin\{2,8,10\}$ would force a root at the same level, contradicting the classification.  Equivalently, repeated reverse jumps would give an impossible infinite descent.
\end{proof}

\subsection{The transition at \texorpdfstring{$\lambda=2$}{lambda=2}}

The matrix of $Q_\lambda$ has eigenvalue
\[
 1+\frac{3\lambda}{2}
\]
in the direction of $(1,1,1,1)$ and eigenvalue
\[
 1-\frac{\lambda}{2}
\]
with multiplicity three. Hence the quadratic form is positive definite for $\lambda<2$, has rank one at $\lambda=2$, and is indefinite for $\lambda>2$, while its value remains positive when all four coordinates are positive.

At the transition,
\begin{equation}\label{eq:lambda2-square}
 Q_2(x,y,z,w)=(x+y+z+w)^2.
\end{equation}
The root condition becomes simply $w\le x+y+z$.

\begin{proposition}\label{prop:lambda23}
For $\lambda=2$, there are exactly $25$ sorted roots, and positive integral solutions exist precisely for
\[
 k\in\{1,2,3,4,5,6,8,9,12,16\}.
\]
For $\lambda=3$, there are exactly $20$ sorted roots, and positive integral solutions exist precisely for
\[
 k\in\{1,2,4,6,7,8,9,11,16,17,22\}.
\]
The complete root lists are given in \ref{app:Qlambda-roots}.
\end{proposition}

\begin{proof}
The bounds of Theorem~\ref{thm:root-strip} give $99$ candidate quadruples $(k,x,y,z)$ for $\lambda=2$ and $216$ for $\lambda=3$. Applying the exact discriminant test yields respectively the $25$ and $20$ roots listed in \ref{app:Qlambda-roots}. Substitution verifies every displayed root, while every other positive solution descends to one of the roots at the same value of $k$.
\end{proof}

There is also a uniform root family.  Setting $(x,y,z,w)=(1,1,1,t)$ in~\eqref{eq:Qlambda-level} gives
\[
 t^2+(3\lambda-k)t+3(\lambda+1)=0,
\]
so
\begin{equation}\label{eq:111t-k}
 k=t+3\lambda+\frac{3(\lambda+1)}{t}.
\end{equation}
Thus $t\mid3(\lambda+1)$, and this solution is a root exactly when
\[
 t^2\le3(\lambda+1).
\]
The admissible values of $k$ for $\lambda=2,3$ therefore differ from $\{2,8,10\}$, the set for $\lambda=1$.

\section{Conclusion and further directions}\label{sec:conclusion}

The reciprocal substitutions describe how the polynomial state changes
under coordinate jumps. The divisibility conditions determine which
integer points realize those changes. Under the stated coprimality
conditions, the six binomial exponents evolve by Vieta mutations. We
classify seven fixed exponent states and twenty-two exponent cycles of
period two. At the chosen representatives of these cycles, eighty-four
sign assignments return after six jumps, and seventy-four admit positive
integral realizations with pairwise-coprime coordinates. The ten excluded cases show why a return
of the signs does not by itself ensure positive coordinates.

Height descent applies to the graph of polynomial states and integer
solutions without assuming that the polynomial data are periodic.
A local root need not remain a root after height-preserving jumps, and
a root stable under such jumps need not have minimum height in its
connected component. For the four-variable family with pairwise cross
terms, the finite bounds give complete root classifications at
$\lambda=1,2,3$ and an exact enumeration procedure for every fixed
positive integer $\lambda$.

Proposition~\ref{prop:cluster-periodic} shows that every nonfixed
periodic exponent state whose intermediate exponents remain positive
must satisfy $ade\ne bcf$. Open questions
include higher even exponent periods under this condition, normalized
expressions with negative powers when exponent differences become
negative, positivity under arbitrarily many six-jump returns, and uniform
descriptions of the roots and connected components for general $\lambda$.
Quadratic families with other cross terms also merit study. Further
comparisons with cluster and Laurent transformations should distinguish
the polynomial identities from the divisibility conditions required for
integer jumps.

\section*{Statements and Declarations}

\paragraph{Funding.}
The authors did not receive support from any organization for the submitted work.

\paragraph{Competing interests.}
The authors declare no competing interests that are relevant to the content of this article.

\paragraph{Data availability.}
All mathematical data supporting the results are provided in the main text
and appendices. The main text contains the fixed exponent states,
construction formulas, proofs ruling out positive pairwise-coprime
realizations in ten cases, and root bounds.
\ref{app:cores} gives the unit-fraction solutions and finite checks for the
period classification; \ref{app:cycles} lists the exponent cycles and
returning sign patterns; \ref{app:positive-certificates} provides starting
data for positive integral realizations; and \ref{app:Qlambda-roots}
lists the four-variable roots for $\lambda=1,2,3$.

\paragraph{Use of generative AI.}
The authors used OpenAI ChatGPT to assist with exploratory calculations, checks, and exposition.  The authors reviewed its outputs and take full responsibility for the content of the paper.

\printbibliography

\clearpage
\appendix
\renewcommand{\thesection}{Appendix \Alph{section}}

\clearpage
\section{Unit-fraction solutions for period two}\label{app:cores}

Theorem~\ref{thm:period2-param} reduces period two to
\[
 \left(\frac1m+\frac1U\right)
 \left(\frac1n+\frac1V\right)=G,
 \qquad G=gh\in\{1,2,3,4\},
\]
with $\gcd(m,n)=1$, $n\mid m+U$, and $m\mid n+V$.
The following bounds show that the enumeration is finite and complete. Put
$\alpha:=\min(m,U)$ and $\beta:=\min(n,V)$, and denote the larger entries
in the two pairs by $M$ and $N$. The reciprocal sums give
$\alpha\beta\le4/G$. Also $G\alpha\beta>1$, since equality or a smaller
value is incompatible with
$(\alpha+M)(\beta+N)=G\alpha\beta MN$.
With $a:=G\alpha\beta-1>0$, this identity is equivalent to
\[
 (aM-\alpha)(aN-\beta)=G\alpha^2\beta^2.
\]
Both factors are positive. Hence, for each of the finitely many pairs
$\alpha,\beta$ with $1<G\alpha\beta\le4$, enumerate the positive divisors
$d$ of $G\alpha^2\beta^2$, set
\[
 M=\frac{d+\alpha}{a},\qquad
 N=\frac{G\alpha^2\beta^2/d+\beta}{a},
\]
and retain the integral values with $M\ge\alpha$, $N\ge\beta$.
Assign the two possible orders independently in each pair and impose the
original gcd and divisibility conditions, removing duplicates. This gives
exactly the following ordered quadruples.

\begin{center}
\small
\begin{tabular}{c p{0.83\textwidth}}
\toprule
$G$ & $(m,n,U,V)$\\
\midrule
1 & $(1,2,3,4)$, $(1,2,5,3)$, $(1,3,2,3)$, $(1,3,5,2)$, $(1,4,3,2)$, $(2,1,3,5)$, $(2,1,4,3)$, $(2,3,1,3)$, $(2,3,4,1)$, $(2,5,3,1)$, $(3,1,2,5)$, $(3,1,3,2)$, $(3,2,1,4)$, $(3,2,3,1)$, $(3,4,1,2)$, $(3,5,2,1)$, $(4,1,2,3)$, $(4,3,2,1)$, $(5,2,1,3)$, $(5,3,1,2)$ \\
2 & $(1,1,2,3)$, $(1,1,3,2)$, $(1,2,1,2)$, $(1,2,3,1)$, $(1,3,2,1)$, $(2,1,1,3)$, $(2,1,2,1)$, $(2,3,1,1)$, $(3,1,1,2)$, $(3,2,1,1)$ \\
3 & $(1,1,1,2)$, $(1,1,2,1)$, $(1,2,1,1)$, $(2,1,1,1)$ \\
4 & $(1,1,1,1)$ \\
\bottomrule
\end{tabular}
\end{center}

Thus there are $20,10,4,1$ ordered quadruples for $G=1,2,3,4$, respectively.

\paragraph{Finite check for the odd-period recurrence.}
We give the finite check from Lemma~\ref{lem:sqrt-recurrence} that the
recurrence terms have integer squares.
After the case $u_{-1}=u_1=1$ has been excluded, write
$u_{-1}=\sqrt a$, $u_0=\sqrt m$, $u_1=\sqrt b$.
The bounds in the proof give $a,b\in\{1,2,3,4\}$ and
$\max(a,b)\le m\le24$.
The squares of $u_2$ and $u_{-2}$ are
\[
 mb+a-2\sqrt{mab},\qquad ma+b-2\sqrt{mab}.
\]
Their integrality requires $mab$ to be a square. Positivity and
maximality also require $mb>a$, $ma>b$, and both displayed integers to
lie between $1$ and $m$. Among triples with $(a,b)\ne(1,1)$, these
finite tests leave precisely the following possibilities:
\begin{center}
\begin{tabular}{ccc l}
\toprule
$u_{-1}$ & $u_0$ & $u_1$ & Subsequent terms $u_2,u_3,\ldots$\\
\midrule
$\sqrt2$ & $2$ & $\sqrt2$ & $\sqrt2,\ 0$\\
$\sqrt2$ & $3$ & $\sqrt2$ & $2\sqrt2,\ 1,\ \sqrt2,\ -\sqrt2$\\
$\sqrt3$ & $2$ & $\sqrt3$ & $\sqrt3,\ 1,\ 0$\\
$2$ & $2$ & $2$ & $2,\ 2,\ldots$\\
\bottomrule
\end{tabular}
\end{center}
Only the constant sequence satisfies the required positivity at every step.

\clearpage
\section{The twenty-two exponent cycles and returning sign patterns}\label{app:cycles}

The exponent vectors are listed in the order $p=(a,b,c,d,e,f)$. The sign patterns refer to the representative $p$ in the first exponent column. A sign pattern records $(\alpha,\beta,\gamma)$; for example \texttt{+-+} means $(+1,-1,+1)$. The final column lists precisely the returning sign patterns that do \emph{not} admit a positive pairwise-coprime six-jump realization.

\footnotesize
\begin{longtable}{r p{0.21\textwidth} p{0.21\textwidth} p{0.28\textwidth} p{0.12\textwidth}}
\toprule
No. & $p$ & $T(p)$ & returning signs & no positive realization\\
\midrule
\endfirsthead
\toprule
No. & $p$ & $T(p)$ & returning signs & no positive realization\\
\midrule
\endhead
1 & $(1,1,3,1,6,3)$ & $(2,1,3,2,3,3)$ & \texttt{+++}, \texttt{--+} & -- \\
2 & $(1,1,4,1,6,4)$ & $(3,2,2,3,2,2)$ & \texttt{+++}, \texttt{+-+}, \texttt{-++}, \texttt{--+} & -- \\
3 & $(1,1,6,2,4,3)$ & $(2,3,2,4,2,1)$ & \texttt{+++}, \texttt{++-}, \texttt{-++}, \texttt{-+-} & -- \\
4 & $(1,1,6,3,3,2)$ & $(1,2,3,3,3,1)$ & \texttt{+++}, \texttt{--+} & -- \\
5 & $(1,2,2,1,6,2)$ & $(3,1,4,3,2,4)$ & \texttt{+++}, \texttt{+-+}, \texttt{-++}, \texttt{--+} & \texttt{+-+} \\
6 & $(1,2,3,1,5,3)$ & $(5,3,2,5,1,2)$ & \texttt{+++}, \texttt{+-+} & -- \\
7 & $(1,2,4,2,3,2)$ & $(3,4,2,6,1,1)$ & \texttt{+++}, \texttt{++-}, \texttt{+-+}, \texttt{+--}, \texttt{-++}, \texttt{-+-}, \texttt{--+}, \texttt{---} & -- \\
8 & $(1,3,2,1,5,2)$ & $(5,2,3,5,1,3)$ & \texttt{+++}, \texttt{--+} & -- \\
9 & $(1,3,2,2,4,1)$ & $(2,1,6,4,2,3)$ & \texttt{+++}, \texttt{++-}, \texttt{-++}, \texttt{-+-} & -- \\
10 & $(1,3,3,3,2,1)$ & $(2,3,3,6,1,1)$ & \texttt{+++}, \texttt{--+} & -- \\
11 & $(1,4,2,2,3,1)$ & $(3,2,4,6,1,2)$ & \texttt{+++}, \texttt{++-}, \texttt{+-+}, \texttt{+--}, \texttt{-++}, \texttt{-+-}, \texttt{--+}, \texttt{---} & \texttt{+-+}, \texttt{+--}, \texttt{---} \\
12 & $(2,1,3,1,4,6)$ & $(4,3,1,2,2,2)$ & \texttt{+++}, \texttt{+-+}, \texttt{-++}, \texttt{--+} & -- \\
13 & $(2,1,5,2,3,5)$ & $(3,5,1,3,2,1)$ & \texttt{+++}, \texttt{-++} & -- \\
14 & $(2,2,2,1,3,4)$ & $(6,4,1,3,1,2)$ & \texttt{+++}, \texttt{++-}, \texttt{+-+}, \texttt{+--} & -- \\
15 & $(2,2,3,2,2,3)$ & $(4,6,1,4,1,1)$ & \texttt{+++}, \texttt{++-}, \texttt{+-+}, \texttt{+--}, \texttt{-++}, \texttt{-+-}, \texttt{--+}, \texttt{---} & -- \\
16 & $(2,3,1,1,4,2)$ & $(4,1,3,2,2,6)$ & \texttt{+++}, \texttt{+-+}, \texttt{-++}, \texttt{--+} & \texttt{-++} \\
17 & $(2,4,1,1,3,2)$ & $(6,2,2,3,1,4)$ & \texttt{+++}, \texttt{++-}, \texttt{+-+}, \texttt{+--} & \texttt{+--} \\
18 & $(2,5,1,2,3,1)$ & $(3,1,5,3,2,5)$ & \texttt{+++}, \texttt{-++} & \texttt{-++} \\
19 & $(2,6,1,2,2,1)$ & $(4,2,3,4,1,3)$ & \texttt{+++}, \texttt{++-}, \texttt{+-+}, \texttt{+--}, \texttt{-++}, \texttt{-+-}, \texttt{--+}, \texttt{---} & \texttt{+--}, \texttt{-++}, \texttt{-+-} \\
20 & $(3,1,2,1,3,6)$ & $(3,2,1,1,3,3)$ & \texttt{+++}, \texttt{--+} & -- \\
21 & $(3,3,1,1,2,3)$ & $(6,3,1,2,1,3)$ & \texttt{+++}, \texttt{--+} & -- \\
22 & $(3,3,2,3,1,2)$ & $(3,6,1,3,1,1)$ & \texttt{+++}, \texttt{--+} & -- \\
\bottomrule
\end{longtable}
\normalsize

\clearpage
\section{Explicit positive realizations}\label{app:positive-certificates}

At the representatives $p$ chosen in \ref{app:cycles}, the following table assigns to each returning sign pattern with a positive pairwise-coprime realization one of the constructions introduced in Section~\ref{sec:positive-lifts}. Each entry specifies an initial positive integer triple through the formulas for $\mathsf P(u)$ or $\mathsf R(X,Y,z)$. To verify an entry, first check pairwise coprimality of the initial coordinates and then substitute into the six successive quotient equations, checking exact division and positivity at every step. The polynomial state returns to its initial value, although the initial and final coordinate triples need not coincide. The six constructions cover $67$ cases:
\[
 \mathsf P(2)\ (45),\quad
 \mathsf P(3)\ (4),\quad
 \mathsf R(2,1,1)\ (14),\quad
 \mathsf R(4,3,1)\ (2),\quad
 \mathsf R(5,2,1)\ (1),\quad
 \mathsf R(13,1,6)\ (1).
\]
The initial triples and parameters for the backward construction in the seven remaining cases are listed separately in Table~\ref{tab:residual-seven}.

\footnotesize
\begin{longtable}{r p{0.85\textwidth}}
\toprule
cycle & sign pattern : construction\\
\midrule
\endfirsthead
\toprule
cycle & sign pattern : construction\\
\midrule
\endhead
1 & \texttt{+++}: \(\mathsf P(2)\); \texttt{--+}: \(\mathsf P(2)\) \\
2 & \texttt{+++}: \(\mathsf P(2)\); \texttt{+-+}: \(\mathsf R(13,16,17)\); \texttt{-++}: \(\mathsf P(2)\); \texttt{--+}: \(\mathsf R(2,1,1)\) \\
3 & \texttt{+++}: \(\mathsf P(2)\); \texttt{++-}: \(\mathsf P(2)\); \texttt{-++}: \(\mathsf P(2)\); \texttt{-+-}: \(\mathsf P(2)\) \\
4 & \texttt{+++}: \(\mathsf P(2)\); \texttt{--+}: \(\mathsf R(2,1,1)\) \\
5 & \texttt{+++}: \(\mathsf P(2)\); \texttt{-++}: \(\mathsf P(2)\); \texttt{--+}: \(\mathsf R(2,1,1)\) \\
6 & \texttt{+++}: \(\mathsf P(2)\); \texttt{+-+}: \(\mathsf R(11,16,3)\) \\
7 & \texttt{+++}: \(\mathsf P(2)\); \texttt{++-}: \(\mathsf P(2)\); \texttt{+-+}: \(\mathsf R(13,35,4)\); \texttt{+--}: \(\mathsf R(7,4,3)\); \texttt{-++}: \(\mathsf P(2)\); \texttt{-+-}: \(\mathsf P(2)\); \texttt{--+}: \(\mathsf R(2,1,1)\); \texttt{---}: \(\mathsf R(4,3,1)\) \\
8 & \texttt{+++}: \(\mathsf P(2)\); \texttt{--+}: \(\mathsf R(2,1,1)\) \\
9 & \texttt{+++}: \(\mathsf P(2)\); \texttt{++-}: \(\mathsf P(2)\); \texttt{-++}: \(\mathsf P(2)\); \texttt{-+-}: \(\mathsf P(2)\) \\
10 & \texttt{+++}: \(\mathsf P(2)\); \texttt{--+}: \(\mathsf P(2)\) \\
11 & \texttt{+++}: \(\mathsf P(2)\); \texttt{++-}: \(\mathsf P(2)\); \texttt{-++}: \(\mathsf P(2)\); \texttt{-+-}: \(\mathsf P(2)\); \texttt{--+}: \(\mathsf R(2,1,1)\) \\
12 & \texttt{+++}: \(\mathsf P(2)\); \texttt{+-+}: \(\mathsf R(2,1,1)\); \texttt{-++}: \(\mathsf R(9,17,10)\); \texttt{--+}: \(\mathsf P(3)\) \\
13 & \texttt{+++}: \(\mathsf P(2)\); \texttt{-++}: \(\mathsf R(10,17,39)\) \\
14 & \texttt{+++}: \(\mathsf P(2)\); \texttt{++-}: \(\mathsf P(2)\); \texttt{+-+}: \(\mathsf R(2,1,1)\); \texttt{+--}: \(\mathsf R(5,2,1)\) \\
15 & \texttt{+++}: \(\mathsf P(2)\); \texttt{++-}: \(\mathsf P(2)\); \texttt{+-+}: \(\mathsf R(2,1,1)\); \texttt{+--}: \(\mathsf R(13,1,6)\); \texttt{-++}: \(\mathsf R(35,13,12)\); \texttt{-+-}: \(\mathsf R(4,3,1)\); \texttt{--+}: \(\mathsf P(3)\); \texttt{---}: \(\mathsf P(3)\) \\
16 & \texttt{+++}: \(\mathsf P(2)\); \texttt{+-+}: \(\mathsf R(2,1,1)\); \texttt{--+}: \(\mathsf P(2)\) \\
17 & \texttt{+++}: \(\mathsf P(2)\); \texttt{++-}: \(\mathsf P(2)\); \texttt{+-+}: \(\mathsf R(2,1,1)\) \\
18 & \texttt{+++}: \(\mathsf P(2)\) \\
19 & \texttt{+++}: \(\mathsf P(2)\); \texttt{++-}: \(\mathsf P(2)\); \texttt{+-+}: \(\mathsf R(2,1,1)\); \texttt{--+}: \(\mathsf P(2)\); \texttt{---}: \(\mathsf P(2)\) \\
20 & \texttt{+++}: \(\mathsf P(2)\); \texttt{--+}: \(\mathsf R(2,1,1)\) \\
21 & \texttt{+++}: \(\mathsf P(2)\); \texttt{--+}: \(\mathsf P(3)\) \\
22 & \texttt{+++}: \(\mathsf P(2)\); \texttt{--+}: \(\mathsf R(2,1,1)\) \\
\bottomrule
\end{longtable}
\normalsize

\begin{table}[ht]
\centering
\small
\caption{Initial positive triples for the seven remaining cases. The column $(X,Y,z)$ records the parameters of the construction $\mathsf R(X,Y,z)$; the last column gives the resulting initial triple.}
\label{tab:residual-seven}
\begin{tabular}{r c c c c}
\toprule
cycle & $p$ & signs & $(X,Y,z)$ & initial triple\\
\midrule
2 & $(1,1,4,1,6,4)$ & \texttt{+-+} & $(13,16,17)$ & $(115,1478,17)$ \\
6 & $(1,2,3,1,5,3)$ & \texttt{+-+} & $(11,16,3)$ & $(7,68,3)$ \\
7 & $(1,2,4,2,3,2)$ & \texttt{+-+} & $(13,35,4)$ & $(55,699,4)$ \\
7 & $(1,2,4,2,3,2)$ & \texttt{+--} & $(7,4,3)$ & $(61,418,3)$ \\
12 & $(2,1,3,1,4,6)$ & \texttt{-++} & $(9,17,10)$ & $(151,37,10)$ \\
13 & $(2,1,5,2,3,5)$ & \texttt{-++} & $(10,17,39)$ & $(572641,2393,39)$ \\
15 & $(2,2,3,2,2,3)$ & \texttt{-++} & $(35,13,12)$ & $(82859,1703,12)$ \\
\bottomrule
\end{tabular}
\end{table}

\clearpage
\section{Root tables for \texorpdfstring{$Q_\lambda=kxyzw$}{Q-lambda levels}}\label{app:Qlambda-roots}

The table gives the complete sorted root lists for $\lambda=1,2,3$.  They are obtained by the finite procedure following Theorem~\ref{thm:root-strip}: enumerate
\[
 \lambda+3\le n=kxy\le6\lambda+4,
\]
then $x\le y$ with $xy\mid n$, impose~\eqref{eq:root-lower-z}--\eqref{eq:root-upper-z}, and apply the exact discriminant test for $w$.

\small
\begin{longtable}{c r p{0.74\textwidth}}
\toprule
$\lambda$ & $k$ & sorted roots $(x,y,z,w)$\\
\midrule
\endfirsthead
\toprule
$\lambda$ & $k$ & sorted roots $(x,y,z,w)$\\
\midrule
\endhead
1 & 2 & \((1,2,7,7)\), \((2,2,2,4)\) \\
1 & 8 & \((1,1,1,2)\) \\
1 & 10 & \((1,1,1,1)\) \\
2 & 1 & \((1,5,24,30)\), \((1,6,14,21)\), \((1,8,9,18)\), \((1,9,10,10)\), \((2,3,10,15)\), \((2,4,6,12)\), \((2,5,5,8)\), \((3,3,6,6)\), \((4,4,4,4)\) \\
2 & 2 & \((1,3,8,12)\), \((1,4,5,10)\), \((2,2,4,8)\), \((2,3,3,4)\) \\
2 & 3 & \((1,2,6,9)\), \((1,3,4,4)\), \((2,2,2,6)\) \\
2 & 4 & \((1,2,3,6)\), \((2,2,2,2)\) \\
2 & 5 & \((1,1,8,10)\), \((1,2,2,5)\) \\
2 & 6 & \((1,1,4,6)\) \\
2 & 8 & \((1,1,2,4)\) \\
2 & 9 & \((1,1,2,2)\) \\
2 & 12 & \((1,1,1,3)\) \\
2 & 16 & \((1,1,1,1)\) \\
3 & 1 & \((1,6,41,49)\), \((1,7,23,32)\), \((1,8,17,27)\), \((1,9,14,25)\), \((1,11,11,24)\), \((1,11,12,17)\), \((2,3,31,31)\), \((3,3,9,12)\), \((3,5,5,8)\), \((4,4,4,12)\) \\
3 & 2 & \((2,2,8,10)\) \\
3 & 4 & \((2,2,2,6)\) \\
3 & 6 & \((1,1,12,13)\) \\
3 & 7 & \((1,1,6,7)\) \\
3 & 8 & \((1,1,4,5)\) \\
3 & 9 & \((1,1,3,4)\) \\
3 & 11 & \((1,1,2,3)\) \\
3 & 16 & \((1,1,1,3)\) \\
3 & 17 & \((1,1,1,2)\) \\
3 & 22 & \((1,1,1,1)\) \\
\bottomrule
\end{longtable}
\normalsize

\end{document}